\documentclass[reqno, 10pt]{article}
	\usepackage{fullpage}
	\usepackage{amsmath}
	\usepackage{amsfonts}
	\usepackage{amssymb}
	\usepackage{amsthm}
	\usepackage{enumerate}
	\usepackage{color}
	\usepackage{bbm}
	\usepackage{cleveref}

\newtheorem{theorem}{Theorem}[section]

\newtheorem{lemma}[theorem]{Lemma}
\newtheorem{proposition}[theorem]{Proposition}

\theoremstyle{definition}
\newtheorem{definition}[theorem]{Definition}

\theoremstyle{remark}
\newtheorem{remark}{Remark}[section]
\newtheorem{nremark}{Notational Remark}[section]

\numberwithin{equation}{section} 

\newcommand{\e}[1]{\mathrm{e}^{#1}}
\newcommand{\R} {\mathbb{R}}
\newcommand{\C} {\mathbb{C}}
\newcommand{\N} {\mathbb{N}}

\newcommand{\dist} {\mathrm{dist}}
\newcommand{\CR} {\C\setminus\R_{+}}

\DeclareMathOperator{\supp}{supp}

\DeclareMathOperator{\re}{\mathrm{Re}}
\DeclareMathOperator{\im}{\mathrm{Im}}

\newcommand{\caD}{{\mathcal D}}

\newcommand{\caH}{{\mathcal H}}

\newcommand{\caV}{{\mathcal V}}

\newcommand{\bsC}{{\boldsymbol C}}

\newcommand{\wt}{\widetilde}
\newcommand{\ol}{\overline}
\newcommand{\wh}{\widehat}
\newcommand{\sa}{\sphericalangle}

\newcommand{\beq}{ \begin{equation} }
\newcommand{\eeq}{ \end{equation} }
\newcommand{\beqs}{	\begin{equation*}	}
\newcommand{\eeqs}{	\end{equation*}	}

\newcommand{\lone}{\mathbbm{1}} 

\newcommand{\dd}{\mathrm{d}}
\newcommand{\ii}{\mathrm{i}}

\newcommand\Var[1]{\mathrm{Var}\left[#1\right]}

\newcommand\absv[1]{\left\vert#1\right\vert}

\newcommand\AND{\quad\text{and}\quad}

\newcommand\ntlim[1]{\lim_{\substack{#1 \\ \sa}}}

\title{Regularity Properties of Free Multiplicative Convolution on the Positive Line}
\author{Hong Chang Ji\footnote{Department of Mathematical Sciences, KAIST, Daejeon, 34141, Korea
		\newline email: \texttt{hcji@kaist.ac.kr}}}

\begin{document}
	
\maketitle

\begin{abstract}
	Given two nondegenerate Borel probability measures $\mu$ and $\nu$ on $\R_{+}=[0,\infty)$, we prove that their free multiplicative convolution $\mu\boxtimes\nu$ has zero singular continuous part and its absolutely continuous part has a density bounded by $x^{-1}$. When $\mu$ and $\nu$ are compactly supported Jacobi measures on $(0,\infty)$ having power law behavior with exponents in $(-1,1)$, we prove that $\mu\boxtimes\nu$ is another Jacobi measure whose density has square root decay at the edges of its support. 
\end{abstract}

\renewcommand{\thefootnote}{\fnsymbol{footnote}} 
\footnotetext{\emph{MSC2010 subject classifications.} 46L54, 30D40, 60B20.}
\renewcommand{\thefootnote}{\arabic{footnote}} 

\section{Introduction}

The notion of free independence, introduced by Voiculescu in \cite{Voiculescu1986}, has been the main object of numerous recent studies, especially after the discovery of its connection to random matrix theory in \cite{Voiculescu1991}. Specifically, free probability provides us a method to calculate the limiting spectral distribution of random matrix ensembles of the form $X+U^{*}YU$ and $X^{1/2}U^{*}YUX^{1/2}$ where $X$ and $Y$ are diagonal matrices having prescribed limiting spectral distributions and $U$ is the Haar unitary matrix. The limiting distributions of such ensembles are given by the free additive and multiplicative convolutions (denoted respectively by $\boxplus$ and $\boxtimes$), which are the distributions of the sum and the product of two free independent random variables with given distributions, and therefore the convolutions themselves have been extensively studied. The purpose of this note is to derive properties of free multiplicative convolution of two probability measures on $[0,\infty)$, namely its absolute continuity and typical behavior of its density.

As shown in \cite{Belinschi2005,Belinschi2006,Belinschi2008}, one of the most interesting features of the free additive and multiplicative convolutions is their absolute continuity except for finitely many point masses, regardless of the measures constituting convolutions unless one of them is degenerate. The following typical example shows how this phenomenon distinguishes the free multiplicative convolution from the classical convolution: consider the probability measure $\mu=\delta_{0}/2+\delta_{2}/2$ and two random variables $X$ and $Y$ both having distribution $\mu$. If $X$ and $Y$ are classically independent, the distribution of $XY$ is $3\delta_{0}/4+\delta_{4}/4$, whereas the free multiplicative convolution $\mu\boxtimes\mu$, the distribution of $X^{1/2}YX^{1/2}$ when $X$ and $Y$ are free, has absolutely continuous part and is explicitly given by
\beq
(\mu\boxtimes\mu)(\dd x)=\frac{1}{2}\delta_{0}(\dd x)+\frac{1}{2\pi}\sqrt{\frac{1}{x(4-x)}}\lone_{(0,4)}(x)\dd x.
\eeq

For generic measures $\mu$ and $\nu$, Voiculescu \cite{Voiculescu1987} derived the equation $S_{\mu\boxtimes\nu}=S_{\mu}S_{\nu}$, where $S_{\mu}$ is defined as
\begin{align}
S_{\mu}(z)&= \frac{(1+z)\chi_{\mu}(z)}{z}, & \chi_{\mu}(\psi_{\mu}(z))&=z, &\psi_{\mu}(z)&=\int_{\R_{+}}\frac{zx}{1-zx}\dd\mu(x).
\end{align}
Since the definition of $S_{\mu}$ involves the inverse mapping of the Cauchy-Stieltjes transforms, calculating $\mu\boxtimes\nu$ amounts to solving the equation satisfied by its Stieltjes transform, and the equation itself is complicated as it also includes the $S$-transforms of $\mu$ and $\nu$. Nevertheless, a complex analytic method that can bypass such complication was first introduced in \cite{Voiculescu1993}, often referred to as analytic subordination. To be precise, given two probability measures $\mu$ and $\nu$ on $\R$, there exist analytic self-maps $\omega_{1}$ and $\omega_{2}$ (referred to as subordination functions) of the upper half plane $\C_{+}$ satisfying $m_{\mu\boxplus\nu}=m_{\mu}\circ\omega_{1}=m_{\nu}\circ\omega_{2}$ in $\C_{+}$, where
\beq\label{eq:Stieltjes}
m_{\mu}(z)= \int_{\R}\frac{1}{x-z}\dd\mu(x),\quad z\in\C_{+}
\eeq
denotes the Stieltjes transform of $\mu$. This result was first proved by Voiculescu for compactly supported measures \cite{Voiculescu1993}, and then by Biane in full generality \cite{Biane1998}. In the same paper Biane also proved the subordination for free multiplicative convolution on $\R_{+}$. (See Theorem~\ref{1.2thm:subor} for its precise statement.) Later, Belinschi and Bercovici proved in \cite{Belinschi-Bercovici2007} that the analytic subordination functions $\omega_{1}$ and $\omega_{2}$ are characterized as solutions of fixed point equations, which is widely used in the present paper. Using the characterization, instead of $S_{\mu\boxtimes\nu}=S_{\mu}S_{\nu}$ we can focus on the equations for $\omega_{1}$ and $\omega_{2}$ that does not involve inverse mappings. The method of subordination has been shown to be the most effective way to prove regularity results for free convolution, for examples, in \cite{Bao-Erdos-Schnelli2018,Belinschi2003,Belinschi2005,Belinschi2006,Belinschi2008}.

In Theorem~\ref{thm:lebesgue}, we follow ideas in \cite{Belinschi2008} to show that $\mu\boxtimes\nu$ has zero singular continuous part. In Theorem~\ref{thm:bound} we use the same ideas as in \cite{Belinschi2014} to prove that the absolutely continuous part of $\mu\boxtimes\nu$ has a continuous density bounded by $1/x$ on $(0,\infty)$, provided that $\mu$ and $\nu$ are nondegenerate and $m_{\mu},m_{\nu}$ are continuous at $0$ and $\infty$. Finally, in Theorem \ref{thm:main}, we follow the line of thoughts in \cite{Bao-Erdos-Schnelli2018} to show that if $\mu$ and $\nu$ are compactly supported Jacobi measures on $(0,\infty)$ having power law behavior with exponents in $(-1,1)$, then $\mu\boxtimes\nu$ is also a Jacobi measure and its density has square root decaying rate at the edges of support. 

All proofs mainly concern the boundary behavior of $M$-functions (see Definition \ref{defn:M} for its definition), another variant of Cauchy-Stieltjes transform, of the measures $\mu$, $\nu$, and $\mu\boxtimes\nu$. To be specific, the function $M_{\mu}$ is defined simply as $\tau\circ \eta_{\mu}\circ \tau$ where $\tau$ is the reciprocal $z\mapsto z^{-1}$ and
\beq
\eta_{\mu}(z)= \frac{\psi_{\mu}(z)}{1+\psi_{\mu}(z)}
\eeq
is the analytic self map of $\CR$ which has been conventionally used in free probability; see \cite{Belinschi-Bercovici2007} for example. Despite of being a simple conjugate of previously known transform, by introducing $M_{\mu}$, we find a particular similarity between the free additive and multiplicative convolutions along the proofs of Theorems \ref{thm:bound} and \ref{thm:main}. For example, $M$-functions enable us to directly connect the main technical result of Belinschi (Lemma \ref{lem:Belinschi} in the present paper), derived along his study of $\mu\boxplus\nu$, to our subject $\mu\boxtimes\nu$ in Section \ref{sec:bound}. Also, we expect the $M$-function to be useful in later researches such as applications in random matrix theory, for it is more directly related to the Stieltjes transform than $\eta$-transform.

Although $M$-functions make the similarity between free additive and multiplicative convolutions more evident, they are fundamentally different and some difficulties arise if one attempts to directly apply arguments for $\mu\boxplus\nu$ in \cite{Bao-Erdos-Schnelli2018,Belinschi2008,Belinschi2014} to $\mu\boxtimes\nu$. An example is the singularity at zero; the equations for the analytic subordination functions $\Omega_{1}$ and $\Omega_{2}$ in Proposition \ref{prop:newsubor} corresponding to $\mu\boxtimes\nu$ become degenerate as $z$ tends to zero. Consequently the behavior of $\mu\boxtimes\nu$ around zero is very different from other positive points (see Theorem~\ref{thm:lebesgue} (ii) and Remark~\ref{rem:sqrt0} for instances). In Theorem~\ref{thm:bound} we deal with zero in a similar way to infinity and prove the bound $1/x$ for the density of $\mu\boxtimes\nu$, in contrast to the constant bound for $\mu\boxplus\nu$ in \cite{Belinschi2014}. On the other hand in Theorem \ref{thm:main}, we simply assume that $\mu$ and $\nu$ are separated from zero so that $\mu\boxtimes\nu$ also is, and this turns out to be necessary for Theorem~\ref{thm:main} to be true at the lower edge. Another difficulty stems from the nonlinearity of analytic subordination in Proposition \ref{prop:newsubor}; we have to closely look at $M_{\mu}/z$ and $\Omega_{j}/z$ while in \cite{Bao-Erdos-Schnelli2018,Belinschi2008} the authors only had to consider $F_{\mu}-z$ and $\omega_{j}-z$, where $\omega_{1}$ and $\omega_{2}$ are analytic subordination functions corresponding to $\mu\boxplus\nu$ and $F_{\mu}=-1/m_{\mu}$. Since $\omega\mapsto\omega/z$ is nonlinear unlike $\omega\mapsto\omega-z$, this makes our analysis more involved, including frequent interchange between $\arg\Omega_{j}$ and $\im\Omega_{j}$.

The paper consists of 5 sections in total. In Section 2, we present preliminary results on general analytic functions on $\C_{+}$ and the analytic subordination in free multiplicative convolution. In Section 3, we formally state our main results, Theorems~\ref{thm:lebesgue}, \ref{thm:bound}, and \ref{thm:main}. The same section also contains the proof of Theorem \ref{thm:lebesgue}. Section 4 is dedicated to the proof of Theorem~\ref{thm:bound}. Finally, Theorem~\ref{thm:main} is proved in Section 5.

\begin{nremark}
	Throughout the paper, we denote the closed positive real axis $[0,\infty)$ by $\R_{+}$ and the open upper half-plane $\{x+\ii y\in\C:y>0\}$ by $\C_{+}$. Unless otherwise indicated, for any subset $A$ of the Riemann sphere $\wh{\C}=\C\cup\{\infty\}$, $\ol{A}$ denotes its closure in $\wh{\C}$. The branch of logarithm is chosen so that $z\mapsto \log z$ maps $\C\setminus(-\infty,0]$ into $\{x+\ii y:x\in\R,y\in(-\pi,\pi)\}$, and the branch of $\arg $ is chosen via $\arg z=\im \log z$ .
\end{nremark}
\begin{nremark}
	For a Borel probability measure $\mu$ on $\R$, we denote its Lebesgue decomposition by $\mu=\mu^{\mathrm{pp}}+\mu^{\mathrm{sc}}+\mu^{\mathrm{ac}}$, where $\mu^{\mathrm{pp}}$, $\mu^{\mathrm{sc}}$, and $\mu^{\mathrm{ac}}$ are the point mass, singular continuous, and absolutely continuous parts of $\mu$, respectively. 
\end{nremark}

\section{Preliminaries}

\begin{definition}\label{defn:M}
	For a probability measure $\mu\neq\delta_{0}$ on $\R_{+}$, we define the functions
	\beq
	M_{\mu}(z)=\frac{1}{\eta_{\mu}(1/z)} \AND
	H_{\mu}(z)= \frac{M_{\mu}(z)}{z},
	\qquad z\in\CR.
	\eeq
\end{definition}

\begin{remark}\label{rem:Mprelim}\
	\begin{itemize}
		\item The maps $M_{\mu}$ and $H_{\mu}$ are analytic in the domain $\CR$.
		\item 
		Note that 
		\beq
		M_{\mu}(z)=\frac{zm_{\mu}(z)}{zm_{\mu}(z)+1}, \quad z\in\CR,
		\eeq
		$M_{\mu}(\C_{+})\subset\C_{+}$, and $M_{\mu}(\ol{z})=\ol{M_{\mu}(z)}$ for $z\in\CR$.
		
		\item  The function $M_{\mu}$ is strictly increasing analytic function on $(-\infty,0)$ as
		\beq
		M_{\mu}'(z)=\left(\int_{\R_{+}}\frac{x}{x-z}\dd\mu(x)\right)^{-2}\left(\int_{\R_{+}}\frac{x}{(x-z)^{2}}\dd\mu(x)\right)>0,
		\quad\forall z\in (-\infty,0).
		\eeq
		Combining this fact with
		\begin{align}\label{eq:ntlimM}
		&\lim_{z\to-\infty,z<0}M_{\mu}(z)=-\infty, &&\lim_{z\to 0^{-}}M_{\mu}(z)=-\frac{\mu(\{0\})}{1-\mu(\{0\})},
		\end{align}
		we find that $M_{\mu}$ restricted to $(-\infty,0)$ has an analytic, strictly increasing inverse on the interval $(-\infty,\frac{-\mu(\{0\})}{1-\mu(\{0\})})$. Using the inverse $M_{\mu}^{-1}$, we can rewrite the defining equation $S_{\mu\boxtimes\nu}=S_{\mu}S_{\nu}$ of the free multiplicative convolution as $zM^{-1}_{\mu\boxtimes\nu}(z)=M^{-1}_{\mu}(z)M^{-1}_{\nu}(z)$.
		
		\item Whenever $\mu$ is nondegenerate, the function $H_{\mu}$ also maps $\C_{+}$ into itself as
		\beq\label{eq:hatIcal}
		\im\left[\frac{1}{z}M_{\mu}(z)\right]
		=\im\frac{m_{\mu}(z)}{zm_{\mu}(z)+1}
		=\frac{\im z}{\absv{zm_{\mu}(z)+1}^{2}} \left[\int_{\R_{+}}\frac{1}{\absv{x-z}^{2}}\dd\mu(x)-\absv{\int_{\R_{+}}\frac{1}{x-z}\dd\mu(x)}^{2}\right].
		\eeq
	\end{itemize}
\end{remark}

\begin{lemma}\label{lem:NPrepM}
	Let $\mu$ be a nondegenerate probability measure on $\R_{+}$ with  
	\beq
	\int_{\R_{+}}x\dd\mu(x)=1\AND \int_{\R_{+}}x^{2}\dd\mu(x)<\infty.
	\eeq
	Then there exists a unique finite Borel measure $\wh{\mu}$ on $\R_{+}$ such that
	\beq\label{eq:NPrepM}
	H_{\mu}(z)=1+\int_{\R_{+}}\frac{1}{x-z}\dd\wh{\mu}(x),\quad
	\wh{\mu}(\R_{+})=\Var{\mu}=\int_{\R_{+}} x^{2}\dd\mu(x)-1.
	\eeq
\end{lemma}
\begin{proof}
	Since $H_{\mu}$ maps $\C_{+}$ into itself, it suffices to check that
	\beq
	\sup_{y\geq 1}\absv{y (H_{\mu}(\ii y)-1)}=\sup_{y\geq 1}\absv{M_{\mu}(\ii y)-\ii y}<\infty
	\eeq
	for the existence and uniqueness of $\wh{\mu}$ (see 1.1 in Chapter 3 of \cite{Akhiezer1965}). Since $M_{\mu}$ is continuous in the imaginary half line $\ii[1,\infty)$, it is enough to show that the limit $\lim_{y\to\infty}(M_{\mu}(\ii y)-\ii y)$ exists, which is calculated as follows:
	\begin{multline}\label{eq:subortotal}
	\lim_{y\to\infty}\ii y\left(\frac{M(\ii y)}{\ii y}-1\right)
	=\lim_{y\to\infty}\frac{\ii y m_{\mu}(\ii y)-\ii y+y^{2}m_{\mu}(\ii y)}{\ii y m_{\mu}(\ii y)+1}
	=-\lim_{y\to\infty}\ii y(\ii y m_{\mu}(\ii y)-\ii y+y^{2}m_{\mu}(\ii y))\\
	=-\lim_{y\to\infty}(\ii y)\left[\left(\ii y m_{\mu}(\ii y)+1\right)-\left(\ii y(\ii y m_{\mu}(\ii y)+1)+1\right)\right]	
	=1+\lim_{y\to\infty}\ii y\int_{\R_{+}}\frac{(\ii y)^{2}+\ii y(x-\ii y)+x(x-\ii y)}{x-\ii y}\dd\mu(x)	\\
	=1+\lim_{y\to\infty}\int_{\R_{+}}\frac{x^{2}(\ii y)}{x-\ii y}\dd\mu(x) 
	=1-\int_{\R_{+}} x^{2}\dd\mu(x)
	=-\Var{\mu},
	\end{multline}
	where we used the fact that
	\beq
	\lim_{y\to\infty}\ii ym_{\mu}(\ii y)=-\mu(\R_{+})=-1,\quad \lim_{y\to\infty}\ii y(\ii y m_{\mu}(\ii y)+1)=-\int_{\R_{+}} x\dd\mu(x)=-1.
	\eeq
	Given the representation, the second identity of \eqref{eq:NPrepM} also directly follows:
	\beq
	\wh{\mu}(\R)=\lim_{y\to\infty}-\ii y\int_{\R}\frac{1}{x-\ii y}\dd\wh{\mu}(x)
	=-\lim_{y\to\infty}(M_{\mu}(\ii y)-\ii y)=\Var{\mu}.
	\eeq
	Since $H_{\mu}$ is analytic and taking real values on $(-\infty,0)$, we have $\supp \wh{\mu}\subset \R_{+}$.
\end{proof}

Along the proof of Theorem \ref{thm:main}, we will repeatedly make use of the following integrals:
\begin{definition}\label{def:I}
	For a Borel probability measure $\mu$ on $\R_{+}$ with mean $1$ and finite variance, we define
	\beq
	I_{\mu}(z)=\int_{\R_{+}}\frac{x}{\absv{x-z}^{2}}\dd\mu(x)\AND
	\wh{I}_{\mu}(z)=\int_{\R_{+}}\frac{1}{\absv{x-z}^{2}}\dd\wh{\mu}(x),
	\eeq
	whenever $z$ is not in the support of the measure in each integral.
\end{definition}

\begin{theorem}[Theorem 3.3 of \cite{Belinschi-Bercovici2007}]\label{1.2thm:subor}
	For two probability measures $\mu$ and $\nu$ on $\R_{+}$ both not $\delta_{0}$, there exist unique analytic functions $\omega_{1},\omega_{2}:\C\setminus\R_{+}\to\C\setminus\R_{+}$ such that
	\begin{itemize}
		\item[(i)] 
		$\lim_{z\to 0^{-}}\omega_{1}(z)=0$ and $\lim_{z\to0^{-}}\omega_{2}(z)=0$,
		
		\item[(ii)] 
		$\omega_{1}$ maps $\C_{+}$ into $\C_{+}$, and for every $z\in\C_{+}$ we have $\omega_{1}(\ol{z})=\ol{\omega_{1}(z)}$ and $\arg\omega_{1}(z)\geq\arg z$. The same statements hold also for $\omega_{2}$. 
		
		\item[(iii)] 
		For any $z\in \C\setminus\R_{+}$, $\eta_{\mu\boxtimes\nu}(z)=\eta_{\mu}(\omega_{1}(z))=\eta_{\nu}(\omega_{2}(z))$.
		
		\item[(iv)] For any $z\in\CR$,
		$\omega_{1}(z)\omega_{2}(z)=z\eta_{\mu\boxtimes\nu}(z)$.
	\end{itemize}
\end{theorem}

\begin{proposition}\label{prop:newsubor}
	Let $\mu$ and $\nu$ be probability measures on $\R_{+}$, both not $\delta_{0}$. Define two analytic functions $\Omega_{1}$ and $\Omega_{2}$ on $\C\setminus \R_{+}$ by
	\beq
	\Omega_{1}(z)= \frac{1}{\omega_{1}(1/z)} \AND
	\Omega_{2}(z)=\frac{1}{\omega_{2}(1/z)},
	\eeq
	where $\omega_{1},\omega_{2}$ are as given in Theorem~\ref{1.2thm:subor}. Then the functions $\Omega_{1}$ and $\Omega_{2}$ map $\C_{+}$ into itself and satisfy the following properties:
	\begin{itemize}
		\item[(i)]
		$\lim_{z\to\infty,z<0}\Omega_{1}(z)=\lim_{z\to\infty,z<0}\Omega_{2}(z)=\infty$ ,	
		\item[(ii)] 
		For any $z\in\C_{+}$, we have $\Omega_{1}(\ol{z})=\ol{\Omega_{1}(z)}$ and $ \arg\Omega_{1}(z)\geq \arg z$, and the same holds for $\Omega_{2}$.
		\item[(iii)] 
		For any $z\in\CR$, $M_{\mu}(\Omega_{1}(z))=M_{\nu}(\Omega_{2}(z))=M_{\mu\boxtimes\nu}(z)$	
		\item[(iv)] 
		For any $z\in\CR$, $\Omega_{1}(z)\Omega_{2}(z)=zM_{\mu\boxtimes\nu}(z).$
	\end{itemize}
\end{proposition}
\begin{proof}
	We prove (i) and (ii) only for $\Omega_{1}$, and the proof is exactly the same for $\Omega_{2}$. First we observe that $\Omega_{1}(z)=1/\omega_{1}(1/z)\in\C_{+}$ whenever $z\in\C_{+}$. For (i) we write
	\beq
	\lim_{z\to\infty,z<0}\Omega_{1}(z)=\lim_{w\to 0,w<0}\Omega_{1}(1/w) =\lim_{w\to 0,w<0}\frac{1}{\omega_{1}(w)}=\infty,
	\eeq
	where the last equality is considered in $\wh{\C}$.
	
	For (ii), from the fact that $\omega_{1}(\bar{z})=\ol{\omega_{1}(z)}$, we have
	\beq
	\arg\Omega_{1}(z)=-\arg\omega_{1}(1/z)=\arg\omega_{1}(1/\bar{z})\geq\arg(1/\bar{z})=z,
	\eeq
	whenever $z\in\C_{+}$.
	
	Identities (iii) and (iv) also follow from those for $\omega_{1}$ and $\omega_{2}$:
	\begin{align}
	M_{\mu}(\Omega_{1}(z))&=\eta_{\mu}\left(\frac{1}{\Omega_{1}(z)}\right)^{-1}=\eta_{\mu}\left(\omega_{1}(1/z)\right)^{-1}=\eta_{\mu\boxtimes\nu}(1/z)^{-1}=M_{\mu\boxtimes\nu}(z),\\
	\Omega_{1}(z)\Omega_{2}(z)&=\frac{1}{\omega_{1}(1/z)\omega_{2}(1/z)}=\left(\frac{1}{z}\eta_{\mu\boxtimes\nu}(1/z)\right)^{-1}=zM_{\mu\boxtimes\nu}(z).
	\end{align}
\end{proof}

\begin{remark}\label{rem:Iim}
	Combining Lemma \ref{lem:NPrepM} and Proposition \ref{prop:newsubor}, we have
	\begin{align}
	\frac{\Omega_{1}(z)}{z}&=m_{\wh{\nu}}(\Omega_{2}(z))+1,&
	\frac{\Omega_{2}(z)}{z}&=m_{\wh{\mu}}(\Omega_{1}(z))+1
	\end{align}
	for $z\in\C_{+}$. In particular, taking imaginary parts of both sides we get
	\begin{align}
	\wh{I}_{\mu}(\Omega_{1}(z))&=\frac{\im(\Omega_{2}(z)/z)}{\im\Omega_{1}(z)}, &
	\wh{I}_{\nu}(\Omega_{2}(z))&=\frac{\im(\Omega_{1}(z)/z)}{\im \Omega_{2}(z)}.
	\end{align}
\end{remark}

By the definition of $\Omega_{1}$ and $\Omega_{2}$, we can simply translate the results of \cite{Belinschi2005,Belinschi2006} concerning boundary behaviors of $\omega_{1}$ and $\omega_{2}$ in terms of $\Omega_{1}$ and $\Omega_{2}$.
\begin{lemma}[Remark 3.3 of \cite{Belinschi2006} and Theorem 1.29 of \cite{Belinschi2005}]\label{lem:suborext} Let $\mu$ and $\nu$ be nondegenerate Borel probability measures on $\R_{+}$. For functions $\Omega_{1}$ and $\Omega_{2}$ defined in Proposition \ref{prop:newsubor}, we have the following:
	\begin{itemize}
		\item[(i)]
		The restrictions $\Omega_{1}\vert_{\C_{+}}$ and $\Omega_{2}\vert_{\C_{+}}$ extend continuously to $\ol{\C_{+}}\setminus\{0\}$, with values in $\ol{\C_{+}}$.
		
		\item[(ii)] Let $x\in\R\setminus\{0\}$ and define
		\begin{align}\label{eq:suborreal}
		\Omega_{1}(x)&=\lim_{z\to x,z\in\C_{+}}\Omega_{1}(z), &
		\Omega_{2}(x)&= \lim_{z\to x,z\in\C_{+}}\Omega_{2}(z).
		\end{align}
		If $\Omega_{1}(x)$ or $\Omega_{2}(x)$ is in $\C_{+}$, then $\Omega_{1}$ and $\Omega_{2}$ extend analytically through $x$ from $\C_{+}$, that is, there exists a neighborhood $U\subset\C$ of $x$ such that $\Omega_{1}\vert_{\C_{+}}$ and $\Omega_{2}\vert_{\C_{+}}$ have analytic extensions to $\C_{+}\cup U$.
	\end{itemize}
\end{lemma}

The rest of this section is dedicated to preliminary results on cluster points and nontangential limits of analytic functions on $\C_{+}$. We will extensively use them in the proofs of Theorems \ref{thm:lebesgue} and \ref{thm:bound}, following \cite{Belinschi2008,Belinschi2014}. We first recall the definition of \emph{nontangential limits}:
\begin{definition}
	Let $f:\C_{+}\to\C$ be a function, $c\in\ol{\R}$, and $\ell\in\wh{\C}$. For $\alpha>0$, define
	\beq
	\Gamma_{\alpha,c}= \{x+\ii y\in\C_{+}:\absv{x-c}<\alpha y\}
	\eeq
	if $c\neq \infty$; otherwise
	\beq
	\Gamma_{\alpha,\infty}= \{x+\ii y\in\C_{+}:\absv{x}<\alpha y\}.
	\eeq
	If we have
	\beq
	\lim_{z\to c,z\in\Gamma_{\alpha,c}}f(z)=\ell,\quad \forall \alpha>0,
	\eeq
	we say that $f$ has nontangential limit $\ell$ at $c$ and write 
	\beq
	\sa\lim_{z\to c}f(z)\equiv\ntlim{z\to c}f(z)=\ell.
	\eeq
	
	If $f$ is defined on $\CR$, we say that $f$ has nontangential limit $\ell$ at $0$ \emph{in $\CR$} if $z\mapsto f(z^{2})$ has nontangential limit $\ell$ at $0$. In this case we write
	\beq
	\ntlim{z\to 0,z\in\CR}f(z)=\ell.
	\eeq
	The nontangential limit at $\infty$ in $\CR$ is defined analogously.
\end{definition}

The result below concerns the relationship between $\mu$ and the nontangential limits of $m_{\mu}$ on $\R$:
\begin{lemma}[Lemma 2.17 of \cite{Belinschi2008}]\label{lem:Stieltjesprop}
	Let $\mu$ be a Borel probability measure on $\R$.
	\begin{itemize}
		\item[(i)] 
		For $\mu^{\mathrm{sc}}$-almost all $x\in\R$, the nontangential limit $\sphericalangle\lim_{z\to x}\im m_{\mu}(z)$ is infinite. 
		
		\item[(ii)] 
		For all $x\in\R$, we have $\mu(\{x\})=-\sa\lim_{z\to x}(z-x)m_{\mu}(z)$.
		
		\item[(iii)] 
		If $f(t)=\frac{\dd \mu^{\mathrm{ac}}(t)}{\dd t}$, we have $\pi f(x)=\sa\lim_{z\to x}\im m_{\mu}(z)$ for almost all $x\in\R$.
	\end{itemize}
\end{lemma}

\begin{lemma}[Theorems 2.5--2.7 of \cite{Belinschi2008}]\label{lem:prelim}\
	\begin{itemize}
		\item[(i)] 
		(Fatou) Let $f:\C_{+}\to \ol{\C_{+}}$ be an analytic function. Then the set of points $x\in\R$ at which the nontangential limit of $f$ fails to exist in $\wh{\C}$ is of Lebesgue measure zero.
		
		\item[(ii)] 
		(Privalov) Let $f:\C_{+}\to\C$ be an analytic function. If there exists a set $A\subset\R$ of nonzero Lebesgue measure such that the nontangential limit of $f$ exists and equals to zero at each point of $A$, then $f$ is constantly zero in $\C_{+}$.
		
		\item[(iii)] 
		(Lindel\"{o}f) Let $f:\C_{+}\to\wh{\C}$ be a meromorphic function that omits at least three points in $\wh{\C}$, and let $x\in\ol{\R}$. If there exists a path $\gamma:[0,1)\to\C_{+}$ such that $\lim_{t\to 1}\gamma(t)=x$ and $\ell=\lim_{t\to 1}f(\gamma(t))$ exists in $\wh{\C}$, then the nontangential limit of $f$ at $x$ exists and equals $\ell$. 
	\end{itemize}
\end{lemma}

\begin{remark}\label{rem:ntlim0}
	As a consequence of Lindel\"{o}f's theorem, we can identify the nontangential limit at zero in $\C_{+}$ with that in $\CR$. For an analytic self-map $f$ of $\CR$ and $\ell\in\wh{\C}$, the following are equivalent;
	\begin{itemize}
		\item[(i)]$\displaystyle\ntlim{z\to0}f(z)=\ell,$
		\item[(ii)]$\displaystyle\ntlim{z\to0,z\in\CR}f(z)=\ell,$
		\item[(iii)]$\displaystyle\lim_{t\to 1^{-}}f(\gamma(t))=\ell \quad \text{ for some path }\gamma:[0,1)\to\C_{+} \quad\text{with }\lim_{t\to 1^{-}}\gamma(t)=0,$
		\item[(iv)]$\displaystyle\lim_{t\to 1^{-}}f(\gamma(t))=\ell \quad \text{ for some path }\gamma:[0,1)\to\CR \quad\text{with }\lim_{t\to 1^{-}}\gamma(t)=0$,
	\end{itemize}
	and the same result holds with zero replaced by $\infty$.
	To prove that (i) implies (ii), we apply Lindel\"{o}f's theorem to the function $z\mapsto f(z^{2})$ and the path $\gamma(t)=\e{\ii\pi/4}\sqrt{1-t}$. Conversely, (ii) implies (i) since for all $\alpha>0$ there exists $\beta>0$ with
	\beq
	\Gamma_{\alpha,\infty}\subset \{z^{2}:z\in\Gamma_{\beta,\infty}\}.
	\eeq
	To deduce (ii) from (iv), we can apply Lindel\"{o}f's theorem to the path $\sqrt{\gamma(t)}$, and we omit the rest of the proof. In particular, applying this result to the function $M_{\mu}$ and paths along the negative real axis, we have
	\begin{align}
	&\ntlim{z\to0,z\in\CR}M_{\mu}(z)=\lim_{z\to0^{-}}M_{\mu}(z)=-\frac{\mu(\{0\})}{1-\mu(\{0\})}=\ntlim{z\to 0}M_{\mu}(z),\\
	&\ntlim{z\to\infty,z\in\CR}M_{\mu}(z)=\lim_{z\to\infty,z<0}M_{\mu}(z)=\infty=\ntlim{z\to\infty}M_{\mu}(z).
	\end{align}
\end{remark}

\begin{lemma}[Lemma 2.13 of \cite{Belinschi2008}]\label{lem:JuliaCara}
	Let $F:\C_{+}\to\C_{+}$ be analytic and let $a\in\R$. If
	\beq
	\ntlim{z\to a}F(z)=c\in\R,
	\eeq
	then
	\beq
	\ntlim{z\to a}\frac{F(z)-c}{z-a}=\liminf_{z\to a}\frac{\im F(z)}{\im z},
	\eeq
	where the equality is considered in $\wh{\C}$. Conversely, if 
	\beq
	\liminf_{z\to a}\frac{\im F(z)}{\im z}<\infty,
	\eeq
	then $\sa\lim_{z\to a}F(z)$ exists and belongs to $\R\cup\{\infty\}$. Moreover, if $F$ is not a constant, then we have $\liminf_{z\to a}\im F(z)/\im z>0$.
\end{lemma}

\begin{definition}
	For a function $f:\C_{+}\to\C_{+}$ and $x\in\ol{\R}$, we define the cluster set of $f$ at $x$ as
	\beq
	\bsC(f,x)= \{w\in\wh{\C}:\exists\{z_{n}\}_{n\in\N}\subset\C_{+} \,\text{ such that }\, z_{n}\to x,\,f(z_{n})\to w\}.
	\eeq
\end{definition}
An important property of $\bsC(f,x)$ is that if $f$ is continuous on $\C_{+}$, then $\bsC(f,x)$ is connected in $\wh{\C}$ (see Lemma 1.6 of \cite{Belinschi2006}).
\begin{lemma}[Lemma 2.18 of \cite{Belinschi2008}]\label{lem:Belinschi}
	Let $f$ be a nonconstant analytic self map of $\C_{+}$. Assume that $x\in\ol{\R}$ is so that $\bsC(f,x)\cap\R$ is infinite. Then there exists a non-empty interval $J\subset\bsC(f,x)\cap\R$ such that for any $c\in J$, we can find a sequence $\{z_{n}^{(c)}\}\subset\C_{+}$ satisfying $\lim_{n\to\infty}z_{n}^{(c)}= x$, $f(z_{n}^{(c)})=c+\ii y_{n}$ for all $n\geq 1$, and $\lim_{n\to\infty}y_{n}=0$.
\end{lemma}

\section{Main results}
We now present our main results. The proof of Theorem \ref{thm:lebesgue} is given here  since it is relatively shorter than that of others.
\begin{theorem}\label{thm:lebesgue}
	Let $\mu$ and $\nu$ be nondegenerate Borel probability measures on $\R_{+}$. Then the following hold true:
	\begin{itemize}
		\item[(i)] 
		$\mu\boxtimes\nu$ has an atom at $c\in(0,\infty)$ if and only if there exists $u,v\in(0,\infty)$ with $uv=c$ and $\mu(\{u\})+\nu(\{\nu\})>1$. In this case, $(\mu\boxtimes\nu)(\{c\})=\mu(\{u\})+\nu(\{v\})-1$.
		
		\item[(ii)] 
		$(\mu\boxtimes\nu)(\{0\})=\max(\mu(\{0\}),\nu(\{0\}))$.
		
		\item[(iii)] 
		The singular continuous part $(\mu\boxtimes\nu)^{\mathrm{sc}}$ is zero.
		
		\item[(iv)]
		Let $f:\R_{+}\to\R_{+}$ be the function defined by
		\beq
		f(x)=\begin{cases}
			\dfrac{1}{\pi}\ntlim{z\to x}\limits\im m_{\mu\boxtimes\nu}(z)& \text{ if the limit exists},\\
			0 & \text{ otherwise}.
		\end{cases}
		\eeq
		Then the function $f$ is a density of $(\mu\boxtimes\nu)^{\mathrm{ac}}$ and there exists a closed set $E\subset\R_{+}$ of Lebesgue measure zero such that $f$ is analytic on $\R_{+}\setminus E$.
	\end{itemize}
\end{theorem}
\begin{proof}
	The first two and the last assertions were proved respectively in \cite{Belinschi2006} and \cite{Belinschi2005}, and are included in this theorem for the sake of completeness. For notational simplicity, we denote $\rho=\mu\boxtimes\nu$ below.	
	
	To prove (iii), we show that $\rho^{\mathrm{sc}}\neq 0$ leads to a contradiction. Assuming $\rho^{\mathrm{sc}}\neq0$, by Lemma~\ref{lem:Stieltjesprop}, there must be a Borel measurable subset $\caH$ of $(0,\infty)$ satisfying the following:
	\begin{itemize}
		\item $\rho^{\mathrm{sc}}(\caH)>0$ and $\rho^{\mathrm{ac}}(\caH)=0=\rho^{\mathrm{pp}}(\caH)$.
		\item $\caH$ is uncountable.
		\item The nontangential limit $\sa\lim_{z\to x}\im (zm_{\rho}(z))$ is infinite for any $x\in \caH$.
	\end{itemize}
	Note that the first follows from the definition, the second from $\rho^{\mathrm{sc}}\perp\rho^{\mathrm{pp}}$, and the last follows from Lemma~\ref{lem:Stieltjesprop} (i) applied to the measure $x\dd\rho^{\mathrm{sc}}(x)$. At each point $c\in \caH$, we claim the following:
	\begin{itemize}
		\item[(a)] 
		$\Omega_{1}(c)$ and $\Omega_{2}(c)$ are in $\R$, where $\Omega_{1}(c)$ and $\Omega_{2}(c)$ are defined in Lemma \ref{lem:suborext}.
		\item[(b)] 
		$\Omega_{1}(c)\Omega_{2}(c)=c$ and $\mu(\{\Omega_{1}(c)\})+\nu(\{\Omega_{2}(c)\})=1$.
	\end{itemize}
	
	Assuming the claims (a) and (b), we first see how $\caH$ being uncountable leads to contradiction. Part (b) of the claim implies that at least one of $\{\Omega_{1}(c):c\in \caH\}$ and $\{\Omega_{2}(c):c\in \caH\}$ must be uncountable, and hence either $\{\Omega_{1}(c):c\in \caH,\mu(\{\Omega_{1}(c)\})>0\}$ or $\{\Omega_{2}(c):c\in \caH,\nu(\{\Omega_{2}(c)\})>0\}$ is also uncountable for $\mu$ and $\nu$ are nondegenerate. Since a probability measure cannot have uncountably many atoms, we obtain contradiction.
	
	Thus it suffices to prove (a) and (b) for all $c\in \caH$. We fix a point $c\in\caH$ and first prove (a). If $\Omega_{1}(c)$ is in $\C_{+}$, Lemma~\ref{lem:suborext} implies that both of the subordination functions extend analytically through $c$. Then $M_{\rho}=M_{\mu}\circ\Omega_{1}$ also analytically extends through $c$, and hence
	\beq
	\sa\lim_{z\to c}M_{\rho}(z)=\sa\lim_{z\to c}M_{\mu}(\Omega_{1}(z))=M_{\mu}(\Omega_{1}(c))\in\C_{+},
	\eeq contradicting $\sa\lim_{z\to c}\im (zm_{\rho}(z))=\infty$. Thus we have $\Omega_{1}(c),\Omega_{2}(c)\in\ol{\R}$, and it suffices to show that $\Omega_{1}(c),\Omega_{2}(c)\neq\infty$ to prove (a). To this end, suppose $\Omega_{1}(c)=\infty$. Then we get
	\beq
	1=\lim_{y\to 0^{+}}M_{\rho}(c+\ii y){=}\lim_{y\to 0^{+}}M_{\mu}(\Omega_{1}(c+\ii y)){=}\ntlim{z\to\infty}M_{\mu}(z),
	\eeq
	where we used Proposition~\ref{prop:newsubor} (iii) and applied Lemma~\ref{lem:prelim} (iii) to the function $M_{\nu}$ and the path $\Omega_{1}(c+\ii(1-t))$. This lead to a contradiction by Remark~\ref{rem:ntlim0}. Similarly we obtain $\Omega_{2}(c)\neq\infty$, concluding the proof of (a).
	
	Now we turn to the proof of (b). We first observe that 
	\beq
	\ntlim{z\to c}\absv{\frac{1}{zm_{\rho}(z)+1}}\leq \ntlim{z\to c}\absv{\frac{1}{\im(zm_{\rho}(z)+1)}}=0,
	\eeq
	from which we obtain
	\beq
	c=\ntlim{z\to c}zM_{\rho}(z)=\Omega_{1}(c)\Omega_{2}(c).
	\eeq
	Also, we observe that
	\beq\label{eq:ntlim1sing}
	1=\ntlim{z\to c}M_{\rho}(z)=\ntlim{z\to c}M_{\mu}(\Omega_{1}(z))=\ntlim{z\to \Omega_{1}(c)}M_{\mu}(z),
	\eeq
	where we again used Lemma~\ref{lem:prelim} (iii) in the last equality. Note also that \eqref{eq:ntlim1sing} implies $\Omega_{1}(c),\Omega_{2}(c)>0$, since $\Omega_{1}(c)< 0$ would lead to $M_{\mu}(\Omega_{1}(c))< 0$.
	
	On the other hand, for any fixed $x\in(0,\infty)$ satisfying $\sa\lim_{z\to x}M_{\mu}(z)=1$, from Lemma~\ref{lem:Stieltjesprop} (ii) we get
	\beq\label{eq:JuliaCaraArg0}
	\frac{1}{\mu(\{x\})}
	=-x\ntlim{z\to x}\left((z-x)(zm_{\mu}(z)+1)\right)^{-1}
	=x\ntlim{z\to x}\frac{M_{\mu}(z)-1}{z-x}=x\ntlim{z\to x}\frac{\log M_{\mu}(z)}{z-x},
	\eeq
	where we used $\lim_{z\to 1}\frac{\log z}{z-1}=1$ in the last equality. Then Lemma~\ref{lem:JuliaCara} implies
	\beq\label{eq:JuliaCaraArg}
	x\ntlim{z\to x}\frac{\log M_{\mu}(z)}{z-x}=\liminf_{z\to x}\frac{\im \log M_{\mu}(z)}{\im z/\absv{z}}	
	=\liminf_{z\to x}\frac{\arg M_{\mu}(z)}{\sin \arg z}
	=\liminf_{z\to x}\frac{\arg M_{\mu}(z)}{\arg z}.
	\eeq
	Note that all equalities in \eqref{eq:JuliaCaraArg0} and \eqref{eq:JuliaCaraArg} remain true even if $\mu(\{x\})=0$ as long as $\sa\lim_{z\to x}M_{\mu}(z)=1$, in which case all limits should be understood as $\infty$. Choosing $x=\Omega_{1}(c)$, we get
	\begin{multline}\label{eq:suborJuliaCara1}
	\frac{1}{\mu(\{\Omega_{1}(c)\})}-1
	=\liminf_{z\to \Omega_{1}(c)}\frac{\arg M_{\mu}(z)-\arg z}{\arg z}
	\leq\liminf_{y\to 0^{+}}\frac{\arg M_{\mu}(\Omega_{1}(c+\ii y))-\arg\Omega_{1}(c+\ii y)}{\arg \Omega_{1}(c+\ii y)}	\\
	{=}\liminf_{y\to 0^{+}}\frac{\arg \Omega_{2}(c+\ii y)-\arg (c+\ii y)}{\arg \Omega_{1}(c+\ii y)}	\leq \liminf_{y\to 0^{+}}\frac{\arg\Omega_{2}(c+\ii y)}{\arg\Omega_{1}(c+\ii y)},
	\end{multline}
	where we used Proposition~\ref{prop:newsubor} (iv). By symmetry, we have the corresponding inequality for $1/\nu(\{\Omega_{2}(c)\})$;
	\beq\label{eq:suborJuliaCara2}
	\frac{1}{\nu(\{\Omega_{2}(c)\})}-1\leq \liminf_{y\to 0^{+}}\frac{\arg \Omega_{1}(c+\ii y)}{\arg \Omega_{2}(c+\ii y)}.
	\eeq
	Note that $\mu(\{\Omega_{1}(c)\})=0$ implies
	\beq
	\liminf_{y\to 0^{+}}\frac{\arg \Omega_{2}(c+\ii y)}{\arg \Omega_{1}(c+\ii y)}=\infty,
	\eeq
	so that
	\beq
	\frac{1}{\nu(\{\Omega_{2}(c)\})}-1\leq \liminf_{y\to 0^{+}}\frac{\arg \Omega_{1}(c+\ii y)}{\arg \Omega_{2}(c+\ii y)}=0,
	\eeq
	contradicting the assumption that $\nu$ is nondegenerate. Similarly $\nu(\{\Omega_{2}(c)\})=0$ leads to contradiction, and thus $\mu(\{\Omega_{1}(c)\}),\nu(\{\Omega_{2}(c)\})\in(0,1)$. Now multiplying \eqref{eq:suborJuliaCara1} and \eqref{eq:suborJuliaCara2}, it follows that 
	\beq\label{eq:thm1ineq}
	\left(\frac{1}{\mu(\{\Omega_{1}(c)\})}-1\right)\left(\frac{1}{\nu(\{\Omega_{2}(c)\})}-1\right)\leq 1,
	\eeq
	which implies $\mu(\{\Omega_{1}(c)\})+\nu(\{\Omega_{2}(c)\})\geq 1$. 
	As $\rho$ does not have point mass at $c$, Theorem~\ref{thm:lebesgue} (i) gives us $\mu(\{\Omega_{1}(c)\})+\nu(\{\Omega_{2}(c)\})=1$. This concludes the proof of (b), and hence the proof of Theorem \ref{thm:lebesgue}.
\end{proof}

Our second result is about the continuity and boundedness of the density of the absolutely continuous part $(\mu\boxtimes\nu)^{\mathrm{ac}}$, which we will prove in Section \ref{sec:bound}.
\begin{theorem}\label{thm:bound}
	Let $\mu$ and $\nu$ be nondegenerate Borel probability measures on $\R_{+}$ such that
	\begin{align}\label{eq:Mext}
	&\lim_{\substack{z\to\infty,\\z\in\CR}}M_{\mu}(z)=\infty=\lim_{\substack{z\to\infty,\\z\in\CR}}M_{\nu}(z),&
	&\lim_{\substack{z\to0,\\z\in\CR}}M_{\mu}(z)=-\frac{\mu(\{0\})}{1-\mu(\{0\})}, &\lim_{\substack{z\to0,\\z\in\CR}}M_{\nu}(z)=-\frac{\nu(\{0\})}{1-\nu(\{0\})}.
	\end{align}
	Furthermore, assume that $\mu(\{a\})+\nu(\{b\})<1$ for all $a,b\in (0,\infty)$. Then the density $f$ of $(\mu\boxtimes\nu)^{\mathrm{ac}}$ defined in Theorem \ref{thm:lebesgue} (iv) is continuous on $(0,\infty)$ and $\sup_{x>0}xf(x)$ is finite.
\end{theorem}

\begin{remark}\
	\begin{itemize}
		\item
		Any probability measures $\mu$ and $\nu$ that are compactly supported in $(0,\infty)$ trivially satisfy the conditions of Theorem \ref{thm:bound}, and there are non-compactly supported examples. Using the asymptotic expansion of the error function in 7.7.2 and 7.12.1 of \cite{HBMF2010}, we can prove that the probability measure $\mu$ defined by
		\beq
		\mu(\dd x)=\frac{1}{Z}\frac{\e{-x^{2}}}{x}\lone_{(1,\infty)}\dd x,
		\eeq
		where $Z$ denotes the normalization constant, satisfies \eqref{eq:Mext}. On the other hand, for the free Poisson distribution $\pi$ defined by 
		\beq\label{eq:MP}
		\pi(\dd x)=\frac{1}{2\pi}\sqrt{\frac{4-x}{x}}\lone_{(0,4)}\dd x,
		\eeq
		we can derive the exact formula of $M_{\pi}$ that satisfies \eqref{eq:Mext}. Taking $\nu=(\mu+\pi)/2$, we have $\supp\nu=\R_{+}$ and $\nu$ satisfies \eqref{eq:Mext}. 
		
		\item
		Denote by $f_{n}$ the density function of the $n$-fold free convolution $\pi^{\boxtimes n}$ of the free Poisson measure $\pi$ defined in \eqref{eq:MP}. Then it is easy to see that the limit
		\beq
		\lim_{x\to 0^{+}}\frac{f_{n}(x)}{x^{-1+1/(n+1)}}
		\eeq
		exists, showing that the boundedness result in Theorem \ref{thm:bound} is optimal.
	\end{itemize}
\end{remark}

Finally, the last result states that the density of $\mu\boxtimes\nu$ decays as square root around the edges of its support if $\mu$ and $\nu$ are Jacobi measures. Its proof is postponed to Section \ref{sec:main}.
\begin{theorem}\label{thm:main}
	Let $\mu$ and $\nu$ be  two absolutely continuous probability measures on $(0,\infty)$, whose densities $f_{\mu}$ and $f_{\nu}$ satisfy the following properties:
	\begin{itemize}
		\item[(a)]	$\supp (f_{\mu})=[E_{-}^{\mu},E_{+}^{\mu}]$, $\supp (f_{\nu})=[E_{-}^{\nu},E_{+}^{\nu}]$ for some positive numbers $E_{-}^{\mu}$, $E_{+}^{\mu}$, $E_{-}^{\nu}$, $E_{+}^{\nu}$;
		
		\item[(b)] $\int_{\R_{+}}x\dd\mu(x)=1=\int_{\R_{+}}x\dd\nu(x)$;
		
		\item[(c)]	there exist exponents $t_{\pm}^{\mu},t_{\pm}^{\nu}\in(-1,1)$ and a constant $\mathfrak{a}>1$ such that
		\beq
		\frac{1}{\mathfrak{a}}\leq \frac{f_{\mu}(x)}{(x-E_{-}^{\mu})^{t_{-}^{\mu}}(E_{+}^{\mu}-x)^{t_{+}^{\mu}}}\leq \mathfrak{a}, \quad \text{for a.e. }x\in[E_{-}^{\mu},E_{+}^{\mu}],
		\eeq
		and
		\beq
		\frac{1}{\mathfrak{a}}\leq \frac{f_{\nu}(x)}{(x-E_{-}^{\nu})^{t_{-}^{\nu}}(E_{+}^{\nu}-x)^{t_{+}^{\nu}}}\leq \mathfrak{a}, \quad\text{for a.e. }x\in[E_{-}^{\nu},E_{+}^{\nu}].
		\eeq
	\end{itemize}
	Then there exist positive numbers $E_{-},E_{+}$, and $\mathfrak{b}>1$ such that the measure $\mu\boxtimes\nu$ is supported on the interval $[E_{-},E_{+}]$ and the density $f$ defined in Theorem \ref{thm:lebesgue} (iv) satisfies
	\beq
	\frac{1}{\mathfrak{b}}\leq \frac{f(x)}{\sqrt{x-E_{-}}\sqrt{E_{+}-x}}\leq \mathfrak{b},\quad \text{for } x\in [E_{-},E_{+}].
	\eeq
\end{theorem}
\begin{remark}\label{rem:sqrt0}\
	\begin{itemize}
		\item Assumption (b) is included for technical simplicity, but it can be easily dropped. Note that the measure $(\mu\boxtimes\nu)$ satisfies the conclusion of Theorem~\ref{thm:main} if and only if $(\mu\boxtimes\nu)^{(a)}$ does for some $a>0$, where $\mu^{(a)}$ denotes the dilation of $\mu$ by a factor $a>0$. Since $\mu^{(a)}\boxtimes\nu^{(b)}=(\mu\boxtimes\nu)^{(ab)}$ for all $a,b>0$, taking $a=(\int_{\R_{+}}x\dd\mu(x))^{-1}$ and $b=(\int_{\R_{+}}x\dd\nu(x))^{-1}$ we find that assumption (b) can be dropped once we prove Theorem~\ref{thm:main}.
		
		\item The assumption that $E_{-}^{\mu}$ and $E_{-}^{\nu}$ are strictly positive is necessary to guarantee the square root decay at the lower edge. In particular if $E_{-}^{\mu}=0=E_{-}^{\nu}$ and $t_{-}^{\mu},t_{-}^{\nu}\in(-1,0)$, applying Karamata's Abelian-Tauberian theorem (see for example Theorem 1.7.4 of \cite{Bingham-Goldie-Teugels1989}) to $M$-functions and their inverses, we can prove from $M_{\mu}^{-1}(w)M_{\nu}^{-1}(w)=wM_{\mu\boxtimes\nu}^{-1}(w)$ that
		\beq\label{eq:sqrt0no}
		(\mu\boxtimes\nu)(0,x)\sim x^{((t_{-}^{\mu}+1)^{-1}+(t_{-}^{\nu}+1)^{-1}-1)^{-1}},
		\eeq
		which is not compatible with the square root decay. Note also that the exponent in \eqref{eq:sqrt0no} is smaller than $\min(t_{-}^{\mu},t_{-}^{\nu})+1$, hence $(\mu\boxtimes\nu)(0,x)$ is much larger than $\mu(0,x)$ and $\nu(0,x)$ when $x$ is small enough.
	\end{itemize}
\end{remark}

\section{Proof of Theorem \ref{thm:bound}}\label{sec:bound}
In this section, we prove Theorem \ref{thm:bound} in two steps. We first prove in Section \ref{sec:bound1} that the subordination functions are continuous at zero, and then use it to prove Theorem \ref{thm:bound} in Section \ref{sec:bound2}.

\subsection{Behavior of subordination functions at $0$}\label{sec:bound1}

\begin{lemma}\label{lem:suborbdry}
	Let $\Omega_{1}$ and $\Omega_{2}$ be the subordination functions defined in Proposition~\ref{prop:newsubor}. Then $C(\Omega_{1},0)\cap\C_{+}=\emptyset=C(\Omega_{2},0)\cap\C_{+}$.
\end{lemma}
\begin{proof}
	Suppose that $\ell\in C(\Omega_{1},0)\cap\C_{+}$, so that there exists a sequence $\{z_{n}\}_{n\in\N}\subset \C_{+}$ which satisfies $\lim_{n\to\infty}z_{n}=0$ and $\lim_{n\to\infty}\Omega_{1}(z_{n})=\ell$. First, we prove that the sequence $\left\{z_{n}H_{\mu}(\Omega_{1}(z_{n}))\right\}_{n\in\N}$ converges to zero nontangentially in $\CR$. Observe that $H_{\mu}(\Omega_{1}(z_{n}))$ converges to $H_{\mu}(\ell)\in\C_{+}$, so that 
	\beq
	\lim_{n\to\infty}\arg H_{\mu}(\Omega_{1}(z_{n}))=\arg H_{\mu}(\ell)\in(0,\pi) \AND \lim_{n\to\infty} z_{n}H_{\mu}(\Omega_{1}(z_{n}))=0.
	\eeq
	In particular there exists $n_{0}\in\N$ so that $\arg H_{\mu}(\Omega_{1}(z_{n}))>(\arg H_{\mu}(\ell))/2$ for all $n\geq n_{0}$. Furthermore the sequence $\{z_{n}H_{\mu}(\Omega_{1}(z_{n}))\}_{n\in\N}$ is identical to $\{\Omega_{2}(z_{n})\}_{n\in\N}$ by Proposition \ref{prop:newsubor}, so that for all $n\geq n_{0}$ we have
	\beq
	\frac{\arg H_{\mu}(\ell)}{2}<\arg H_{\mu}(\Omega_{1}(z_{n}))<\arg z_{n}+\arg H_{\mu}(\Omega_{1}(z_{n}))=\arg(z_{n}H_{\mu}(\Omega_{1}(z_{n})))=\arg \Omega_{2}(z_{n})<\pi.
	\eeq
	Therefore we conclude that $\{z_{n}H_{\mu}(\Omega_{1}(z_{n}))\}_{n\in\N}$ converges to zero nontangentially in $\CR$ since the sequence is contained in the sector $\{\omega\in\C_{+}:\arg H_{\mu}(\ell)/2<\arg \omega<\pi\}$ for $n\geq n_{0}$.
	
	Then, by Remark \ref{rem:ntlim0}, we obtain
	\beq
	\lim_{n\to\infty}M_{\nu}(z_{n}H_{\mu}(\Omega_{1}(z_{n})))=\lim_{\substack{\omega\to0, \\ \arg H_{\mu}(\ell)/4< \arg \omega< \pi/2}}M_{\nu}(\omega^{2})=\ntlim{z\to0,z\in\CR}M_{\nu}(z){=}-\frac{\nu(\{0\})}{1-\nu(\{0\})}.
	\eeq
	Since $M_{\mu}(\Omega_{1}(z))=M_{\nu}(\Omega_{2}(z))$ and $\Omega_{2}(z)=zH_{\mu}(\Omega_{1}(z))$, we conclude that
	\beq\label{eq:subor0im}
	\ell=\lim_{n\to\infty}\Omega_{1}(z_{n})=\lim_{n\to\infty}\frac{M_{\mu}(\Omega_{1}(z_{n}))}{H_{\mu}(\Omega_{1}(z_{n})}=\lim_{n\to\infty}\frac{M_{\nu}(z_{n}H_{\mu}(\Omega_{1}(z_{n})))}{H_{\mu}(\Omega_{1}(z_{n}))}{=}-\frac{\nu(\{0\})}{H_{\mu}(\ell)(1-\nu(\{0\}))}.
	\eeq
	Now multiplying \eqref{eq:subor0im} by $H_{\mu}(\ell)$, we find that
	\beq
	-\frac{\nu(\{0\})}{1-\nu(\{0\})}=M_{\mu}(\ell)\in\C_{+},
	\eeq
	which is a contradiction.
\end{proof}

\begin{lemma}\label{lem:suborntlim}
	The functions $\Omega_{1}$ and $\Omega_{2}$ have nontangential limits at zero and $\infty$ in $\CR$, given by
	\begin{align}
	\ntlim{z\to 0,z\in\CR}\Omega_{1}(z)&=
	\begin{cases}
	0	&\text{if }\mu(\{0\})\geq\nu(\{0\}), \\
	M_{\mu}^{-1}\left(-\frac{\nu(\{0\})}{1-\nu(\{0\})}\right) & \text{if }\nu(\{0\})>\mu(\{0\}),
	\end{cases}\label{eq:ntlim1}\\
	\ntlim{z\to 0,z\in\CR}\Omega_{2}(z)&=
	\begin{cases}
	0	&\text{if }\nu(\{0\})\geq\mu(\{0\}), \\
	M_{\nu}^{-1}\left(-\frac{\mu(\{0\})}{1-\mu(\{0\})}\right) & \text{if }\mu(\{0\})>\nu(\{0\}),
	\end{cases}\label{eq:ntlim2}\\
	\ntlim{z\to\infty,z\in\CR}\Omega_{1}(z)&=\infty=\ntlim{z\to\infty,z\in\CR}\Omega_{2}(z).
	\end{align}
\end{lemma}
\begin{proof}
	Note that by the second item of Remark \ref{rem:Mprelim}, the quantities on the right-hand sides of \eqref{eq:ntlim1} and \eqref{eq:ntlim2} are well-defined. We prove the lemma only for $\Omega_{1}$, and the proof for $\Omega_{2}$ is completely analogous. Applying Remark \ref{rem:ntlim0} to the function $\Omega_{1}$, it follows that
	\begin{align}
	&\ntlim{z\to 0,z\in\CR}\Omega_{1}(z)=\lim_{z\to 0^{-}}\Omega_{1}(z), &
	&\ntlim{z\in\infty,z\in\CR}\Omega_{1}(z)=\lim_{z\to -\infty,z<0}\Omega_{1}(z),
	\end{align}
	provided that the limits on the right-hand sides exist.
	Thus we can restrict our attention to the limits along the negative real axis. 
	
	For $\lim_{z\to 0^{-}}\Omega_{1}(z)$, Theorem~\ref{thm:lebesgue} (ii) implies
	\beq\label{eq:MR-}
	M_{\rho}(-\infty,0)=(-\infty,-\frac{\rho(\{0\})}{1-\rho(\{0\})})=
	(-\infty,-\frac{\mu(\{0\})}{1-\mu(\{0\})})\cap(-\infty,-\frac{\nu(\{0\})}{1-\nu(\{0\})}),
	\eeq
	so that using Proposition~\ref{prop:newsubor} (iii) we get
	\beq
	\Omega_{1}(z)=M_{\mu}^{-1}\circ M_{\rho}(z)
	\eeq
	on the whole open line $(-\infty,0)$ as the right-hand side of \eqref{eq:MR-} is contained in the analytic domain of $M_{\mu}^{-1}$ in Remark \ref{rem:Mprelim}. 
	Thus we obtain
	\beq
	\lim_{z\to 0^{-}}\Omega_{1}(z)=\lim_{z\to\left(-\frac{\rho(\{0\})}{1-\rho(\{0\})}\right)^{-}}M_{\mu}^{-1}(z)=
	\begin{cases}
		0	&\text{if }\mu(\{0\})\geq\nu(\{0\}), \\
		M_{\mu}^{-1}\left(-\frac{\rho(\{0\})}{1-\rho(\{0\})}\right) & \text{if }\nu(\{0\})>\mu(\{0\}),
	\end{cases}
	\eeq
	and similarly
	\beq
	\lim_{z\to -\infty,z\in\R}\Omega_{1}(z)=\lim_{z\to -\infty,z\in\R}M_{\mu}^{-1}(z)=-\infty,
	\eeq
	concluding the proof.
\end{proof}

\begin{proposition}\label{prop:suborext0}
	The restrictions $\Omega_{1}\vert_{\C_{+}}$ and $\Omega_{2}\vert_{\C_{+}}$ extend continuously to $\ol{\C_{+}}$, and the following hold true:
	\begin{align}
	\lim_{z\to0,z\in\C_{+}}\Omega_{1}(z)&=
	\begin{cases}
	0	&\text{if }\mu(\{0\})\geq\nu(\{0\}), \\
	M_{\mu}^{-1}\left(-\frac{\nu(\{0\})}{1-\nu(\{0\})}\right) & \text{if }\nu(\{0\})>\mu(\{0\}),
	\end{cases}\label{eq:suborext01}\\
	\lim_{z\to0,z\in\C_{+}}\Omega_{2}(z)&=
	\begin{cases}
	0	&\text{if }\nu(\{0\})\geq\mu(\{0\}), \\
	M_{\nu}^{-1}\left(-\frac{\mu(\{0\})}{1-\mu(\{0\})}\right) & \text{if }\mu(\{0\})>\nu(\{0\}),
	\end{cases}\\
	\lim_{z\to\infty,z\in\C_{+}}\Omega_{1}(z)&=\infty=\lim_{z\to\infty,z\in\C_{+}}\Omega_{2}(z).\label{eq:suborextinfi}
	\end{align}
\end{proposition}
\begin{proof}
	We prove the result for $\Omega_{1}$, and the proof for $\Omega_{2}$ is analogous. For simplicity, we denote the right-hand side of \eqref{eq:suborext01} by $\Omega_{1}(0)$. First, we prove that \eqref{eq:suborext01} and \eqref{eq:suborextinfi} imply that $\Omega_{1}\vert_{\C_{+}}$ extends continuously to $\ol{\C_{+}}$. For an arbitrarily chosen $\epsilon>0$, \eqref{eq:suborext01} implies that we can find $\delta>0$ so that $\absv{\Omega_{1}(z)-\Omega_{1}(0)}\leq\epsilon$ holds whenever $z\in\C_{+}$ and $\absv{z}<\delta$. Then by \eqref{eq:suborreal}, the same holds for all $z\in\C_{+}\cup \R$ with $\absv{z}<\delta$. Applying the same reasoning with zero replaced by $\infty$, we conclude
	\begin{align}
	&\lim_{z\to0,z\in\ol{\C_{+}}}\Omega_{1}(z)=\Omega_{1}(0),&
	&\lim_{z\to\infty,z\in\ol{\C_{+}}}\Omega_{1}(z)=\infty.
	\end{align}
	Consequently, if we define $\Omega_{1}(x)$, $\Omega_{1}(0)$, and $\Omega_{1}(\infty)$ as in respectively \eqref{eq:suborreal}, \eqref{eq:suborext01}, and \eqref{eq:suborextinfi}, the function $\Omega_{1}:\ol{\C_{+}}\to\wh{\C}$ is continuous.
	
	Next, we prove that $\bsC(\Omega_{1},0)$ is a singleton. Since $\bsC(\Omega_{1},0)$ is a connected subset of $\ol{\R}$ by Lemma \ref{lem:suborbdry}, $\bsC(\Omega_{1},0)$ either is a singleton or contains an interval. Thus it suffices to prove that the cluster set has measure zero. In order to do so, suppose on the contrary that $\bsC(\Omega_{1},0)\subset\ol{\R}$ has positive measure. Then for any $c\in \bsC(\Omega_{1},0)\setminus\{0,\infty\}$ for which the nontangential limit $\sa\lim_{z\to c}M_{\mu}(z)$ exists and is finite, we use Lemma~\ref{lem:Belinschi} to take a sequence $\{z_{n}^{(c)}\}_{n\in\N}$ such that $z_{n}^{(c)}\to 0$, $\Omega_{1}(z_{n}^{(c)})\to c$, $\Omega_{1}(z_{n}^{(c)})\in c+\ii\R_{+}$. Then by Proposition \ref{prop:newsubor} we have
	\beq
	\lim_{n\to\infty}\Omega_{2}(z_{n}^{(c)})=\lim_{n\to\infty}z_{n}^{(c)}\frac{M_{\mu}(\Omega_{1}(z_{n}^{(c)}))}{\Omega_{1}(z_{n}^{(c)})}=0,
	\eeq
	and hence by \eqref{eq:Mext} we obtain
	\beq
	\ntlim{z\to c}M_{\mu}(z)
	=\lim_{n\to\infty}M_{\mu}(\Omega_{1}(z_{n}^{(c)}))
	=\lim_{n\to\infty}M_{\nu}(\Omega_{2}(z_{n}^{(c)}))
	=\lim_{z\to 0}M_{\nu}(z){=}-\frac{\nu(\{0\})}{1-\nu(\{0\})}.
	\eeq
	Since $c$ was chosen from arbitrary from a set of positive measure, by Lemma~\ref{lem:prelim} (ii) we have a contradiction.
	
	Finally, we prove \eqref{eq:suborext01} and \eqref{eq:suborextinfi}. By Lemma \ref{lem:suborext} the set $\bsC(\Omega_{1},\infty)$ is a singleton, and by above $\bsC(\Omega_{1},0)$ is also a singleton. Combining Remark \ref{rem:ntlim0} and Lemma \ref{lem:suborntlim} we find that $\infty\in \bsC(\Omega_{1},\infty)$ and $\Omega_{1}(0)\in\bsC(\Omega_{1},0)$. Therefore we have $\bsC(\Omega_{1},0)=\{\Omega_{1}(0)\}$ and $\bsC(\Omega_{1},\infty)=\{\infty\}$, which imply \eqref{eq:suborext01} and \eqref{eq:suborextinfi}, respectively.
\end{proof}

\subsection{Proof of Theorem~\ref{thm:bound}}\label{sec:bound2}
We will first prove that if $\mu$ and $\nu$ satisfy the assumptions of Theorem \ref{thm:bound}, then $\absv{M_{\rho}-1}$ is bounded below by a strictly positive constant on $\C_{+}$. Then we shall deduce the boundedness of $xf(x)$ from it and prove the continuity. Suppose on the contrary that there exists a sequence $\{z_{n}\}_{n\in\N}$ in $\C_{+}$ such that $\lim_{n\to\infty}M_{\rho}(z_{n})=1$. Taking a subsequence, we assume that the sequence $\{z_{n}\}$ converges to a point $c\in \ol{\C_{+}}$. Furthermore, since $M_{\rho}$ is an analytic self-map of $\C\setminus\R_{+}$, we may assume that $c\in[0,\infty]$. We will show that $c=0$ or $\infty$ leads to contradiction, and $c\in(0,\infty)$ implies that $\mu(\{v\})+\nu(\{u\})\geq 1$ which also contradicts the assumption.

We first assume $c=0$. By Proposition~\ref{prop:suborext0}, we see that either $\lim_{n\to\infty}\Omega_{1}(z_{n})=0$ or $\lim_{n\to\infty}\Omega_{2}(z_{n})=0$ must hold. If the former were true, using \eqref{eq:Mext} we would get
\beq
\lim_{n\to\infty}M_{\rho}(z_{n})=\lim_{n\to\infty}M_{\mu}(\Omega_{1}(z_{n}))=\lim_{z\to 0,z\in\CR}M_{\mu}(z)=-\frac{\mu(\{0\})}{1-\mu(\{0\})}\neq 1.
\eeq
Similarly $\lim_{n\to\infty}\Omega_{2}(z_{n})=0$ leads to a contradiction. Second, we assume $c=\infty$. Again by Lemmas~\ref{lem:suborext} and \ref{lem:suborntlim}, we have $\lim_{n\to\infty}\Omega_{1}(z_{n})=\infty$ so that \eqref{eq:Mext} implies
\beq
\lim_{n\to\infty}M_{\rho}(z_{n})=\lim_{n\to\infty}M_{\mu}(\Omega_{1}(z_{n}))=\lim_{z\to \infty,z\in\CR}M_{\mu}(z)=\infty.
\eeq

Therefore $c$ must be a positive real number. From Lemma~\ref{lem:suborext}, we already know that
\beq
\Omega_{1}(c)=\lim_{z\to c,z\in\C_{+}}\Omega_{1}(z)\AND \Omega_{2}(c)=\lim_{z\to c,z\in\C_{+}}\Omega_{2}(z)
\eeq
both exist, and we will first prove that they are also positive real numbers. As 
\beq
\arg \Omega_{1}(z_{n})\leq \arg M_{\rho}(z_{n})\to 0, 
\eeq
we have $\Omega_{1}(c)\in [0,\infty]$ and the same argument shows $\Omega_{2}(c)\in [0,\infty]$. If $\Omega_{1}(c)$ were infinity, \eqref{eq:Mext} would imply
\beq
\infty=\lim_{z\to\infty}M_{\mu}(z)
=\lim_{n\to\infty}M_{\mu}(\Omega_{1}(z_{n}))
=\lim_{n\to\infty}M_{\rho}(z_{n})=1,
\eeq
leading to contradiction. Similarly $\Omega_{2}(c)<\infty$. 
Also $\Omega_{1}(c)=0$ would imply
\beq
0=\Omega_{1}(c)\Omega_{2}(c)=cM_{\rho}(c)=c,
\eeq
so that $\Omega_{1}(c)$ and $\Omega_{2}(c)$ are both positive real numbers.

Then Proposition~\ref{prop:newsubor} (iv) together with Lemma~\ref{lem:suborext} implies $M_{\rho}$ also continuously extends to $c$ with value $1$. Thus
\beq
\lim_{y\to 0^{+}}M_{\mu}(\Omega_{1}(c+\ii y))=\lim_{y\to 0^{+}}M_{\rho}(c+\ii y)=1,
\eeq
and Lemma~\ref{lem:prelim} (iii) implies existence of the nontangential limits
\beq
\ntlim{z\to \Omega_{1}(c)}M_{\mu}(z)=\ntlim{z\to\Omega_{2}(c)}M_{\nu}(z)=1.
\eeq
Given the existence of nontangential limits, \eqref{eq:JuliaCaraArg} implies
\beq
\ntlim{z\to \Omega_{1}(c)}\frac{M_{\mu}(z)-1}{z-\Omega_{1}(c)}=\liminf_{z\to\Omega_{1}(c)}\frac{\arg M_{\mu}(z)}{\arg z}
\eeq
and similar equality for $\Omega_{2}$ and $M_{\nu}$. Following the lines below \eqref{eq:JuliaCaraArg}, we again obtain
\beq
\mu(\{\Omega_{1}(c)\})+\nu(\{\Omega_{2}(c)\})\geq 1,
\eeq
contradicting our assumption. Therefore we have proved that $\absv{M_{\rho}-1}$ is bounded below by a strictly positive constant on $\C_{+}$, or equivalently, that the set $\{zm_{\rho}(z)+1:z\in\C_{+}\}=\{(1-M_{\rho}(z))^{-1}:z\in\C_{+}\}$ is bounded. Then the Stieltjes inversion(Lemma \ref{lem:Stieltjesprop} (iii)) directly implies that the density $xf(x)$ is bounded. 

Now we prove the continuity of $f$. By Proposition~\ref{prop:suborext0}, both of $\Omega_{1}$ and $\Omega_{2}$ are continuous at $0$ and either of them must have value $0$ at $0$. Thus Proposition~\ref{prop:newsubor} (iii) and \eqref{eq:Mext} imply that $\absv{\bsC(M_{\rho},0)}=1$. For $c>0$, Lemma~\ref{lem:suborext} and Proposition~\ref{prop:newsubor} (iv) implies that $\absv{\bsC(M_{\rho},c)}=1$. Note that we are allowing the value $\infty$ for $M$, but nevertheless $\R\ni t\mapsto tm_{\mu}(t)$ is a continuous bounded function.

\begin{remark}\label{rem:bound0pf}
	Even if $M_{\mu}$ or $M_{\nu}$ is not continuous at $0$, every argument in the proof remains intact except that we cannot exclude the case $c=0$ and that we cannot guarantee that $M_{\rho}$ is continuous at $0$. Thus, in this case also, we have the continuity and boundedness of $xf(x)$ on $(\epsilon,\infty)$ for all $\epsilon>0$.
\end{remark}

\section{Proof of Theorem \ref{thm:main}}\label{sec:main}

In this section we focus on the proof of Theorem~\ref{thm:main}. The exposition of the proof closely follows that of \cite{Bao-Erdos-Schnelli2018}, and most of the lemmas have their counterparts in terms of free additive convolution in \cite{Bao-Erdos-Schnelli2018}. The proof is divided into two steps; the first ingredient is Proposition~\ref{prop:stab}, which states that $\Omega_{1}$ and $\Omega_{2}$ stay away from $\supp \mu$ and $\supp \nu$, respectively. The lower boundedness of $\dist(\Omega_{1}(\C_{+}),\supp\mu)$ and $\dist(\Omega_{2}(\C_{+}),\supp\nu)$ is often called \emph{stability bound}, as it is directly related to the stability of subordination equations in Proposition~\ref{prop:newsubor}: see \cite{Bao-Erdos-Schnelli2016,Lee-Schnelli2013} for examples in case of free additive convolution.  

Given the stability bound, we then prove in Proposition \ref{prop:edgechar2} that the density of $\mu\boxtimes\nu$ has a single interval support $[E_{-},E_{+}]$, and that the density is strictly positive in the interior $(E_{-},E_{+})$. Along the proof, we will find that the edges $E_{\pm}$ are characterized as the points at which the inverse functions of $\Omega_{1}$ and $\Omega_{2}$ have vanishing derivatives. From this fact, in Proposition \ref{prop:sqrt} we deduce that $\im\Omega_{1}$ and $\im\Omega_{2}$ decay as square root at the edges, which will imply the square root decay of the density of $\mu\boxtimes\nu$.

\begin{nremark}
	As in the previous sections, we denote $\rho=\mu\boxtimes\nu$ for simplicity. We denote the density of $\rho$ given in Theorem \ref{thm:lebesgue} by $f$. Unless otherwise specified, we will always denote $\re z=E$ and $\im z=y$ for $z\in\C_{+}$. Finally, we will denote by $C,c$ positive constants independent of $z$, and their value may vary from line to line.
\end{nremark}

An important consequence of our assumptions in Theorem \ref{thm:main} is that $\mu$ and $\wh{\mu}$ have the same support, which will allow us to effectively deal with $m_{\wh{\mu}}$ instead of the reciprocal $(zm_{\mu}(z)+1)^{-1}$ and simplify the proof.
\begin{lemma}\label{lem:NPrepMsupp}
	Let $\mu$ and $\nu$ satisfy the assumptions of Theorem \ref{thm:main}, and let $\wh{\mu}$ and $\wh{\nu}$ be the Borel measures on $\R$ corresponding to $\mu$ and $\nu$ via Lemma \ref{lem:NPrepM}. Then
	\beq\label{eq:NPrepMsupp}
	\supp\mu=\supp\wh{\mu} \AND \supp\nu=\supp\wh{\nu}.
	\eeq 
\end{lemma}
\begin{proof}
	We first prove that $\supp\wh{\mu}\subset\supp \mu$. Since $\mu$ is supported on a single interval, for any $z\in\C\setminus\supp\mu$ either one of $\re z<E_{-}^{\mu}$ or $\im z>0$ or $\re z>E_{+}^{\mu}$ must hold, so that $zm_{\mu}(z)+1\neq0$. Therefore $z\mapsto (M_{\mu}(z)-z)$ is analytic in $\C\setminus\supp\mu$ with $\im(M_{\mu}(z)-z)=0$ for $z\in\R\setminus\supp \mu$. Thus, whenever $x\in\R\setminus\supp\mu$ we have 
	\beq
	\lim_{y\to 0^{+}}\im m_{\wh{\mu}}(x+\ii y)=0,
	\eeq
	which proves $\supp\wh{\mu}\subset\supp \mu$.
	
	We are now left with the reverse inclusion $\supp\mu\subset\supp\wh{\mu}$. Suppose the contrary, and take a non-empty open interval $I\subset\supp\mu\setminus\supp\wh{\mu}$. By Stieltjes inversion, for almost every $E\in I$, we have
	\beq
	\lim_{y\to 0^{+}}\im m_{\mu}(E+\ii y)=\pi f_{\mu}(E)>0
	\eeq
	On the other hand, $I\subset(\supp\wh{\mu})^{\mathrm{c}}$ implies that $z\mapsto\frac{M_{\mu}(z)}{z}$ extends to $I$ with real-values on $I$, so that $z\mapsto(\frac{M_{\mu}(z)}{z}-1)$ extends to an analytic function defined on $D=\{E+\ii y:E\in I,y\in\R\}$ by Schwarz reflection. Also $M_{\mu}-1$ is not identically zero in $D$, for $M_{\mu}-1\equiv0$ would imply $(zm_{\mu}(z)+1)^{-1}=0$ for all $z\in\C_{+}$ which is impossible. Thus there are only finitely many solutions of the equation $M_{\mu}(z)=1$ and the formula
	\beq
	zm_{\mu}(z)+1=\frac{1}{1-M_{\mu}(z)}
	\eeq
	defines a meromorphic extension of $zm_{\mu}(z)+1$ on $D$, which is real-valued on $I$. Therefore for almost every $E\in I\setminus\{0\}$ for which $1-M_{\mu}(E)\neq0$ we get 
	\beq
	\lim_{y\to 0^{+}}m_{\mu}(E+\ii y)=\frac{M_{\mu}(E)}{E(1-M_{\mu}(E))}\in\R,
	\eeq
	implying $\lim_{y\to 0^{+}}\im m_{\mu}(E+\ii y)=0$. Since $E\in I\subset (E_{-}^{\mu},E_{+}^{\mu})$, this leads to a contradiction as $f_{\mu}(E)>0$ by Lemma \ref{lem:prelim} (iii).
\end{proof}

\subsection{Stability bounds}

The goal of this section is to prove the stability bound, Proposition \ref{prop:stab}. Before its statement and proof, we first prove that the subordination functions $\Omega_{1}$ and $\Omega_{2}$ stay away from $0$ and $\infty$ if their argument $z$ does:
\begin{lemma}\label{lem:suborbd}
	Let $\mu$ and $\nu$ be as in Theorem \ref{thm:main} and let $\caD\subset\ol{\C_{+}}\setminus\{0,\infty\}$ be compact. Then there exists a constant $C>1$ depending on $\mu$, $\nu$, and $\caD$ such that for all $z\in \caD$,
	\beq
	C^{-1}\leq \absv{\Omega_{1}(z)}\leq C,\quad 
	C^{-1}\leq \absv{\Omega_{2}(z)}\leq C.
	\eeq
\end{lemma}
\begin{proof}
	By Lemma~\ref{lem:suborext}, it suffices to prove that for any $c\in\R\setminus\{0\}$, $\Omega_{1}(c)$ and $\Omega_{2}(c)$ cannot be zero or infinity. As in the previous proofs, we prove only for $\Omega_{1}$.
	
	Suppose on the contrary that $\Omega_{1}(c)=0$ for some $c\in\R\setminus\{0\}$. Then 
	\beq
	\lim_{y\to 0^{+}}\frac{\Omega_{2}(c+\ii y)}{c+\ii y}
	=\lim_{y\to 0^{+}}\frac{M_{\mu}(\Omega_{1}(c+\ii y))}{\Omega_{1}(c+\ii y)}
	=\lim_{z\to 0}\frac{M_{\mu}(z)}{z}=\int_{\R_{+}}\frac{1}{x}\dd\mu(x)\in (0,\infty),
	\eeq
	which in turn implies $\Omega_{2}(c)$ cannot be $0$ nor $\infty$. We then have 
	\beq
	\lim_{y\to0^{+}}M_{\nu}(\Omega_{2}(c+\ii y))=\lim_{y\to0^{+}}M_{\rho}(c+\ii y)=\lim_{y\to0^{+}}\frac{\Omega_{1}(c+\ii y)\Omega_{2}(c+\ii y)}{c+\ii y}= 0,
	\eeq
	as $y\to 0^{+}$. This implies, from the definition of $M_{\rho}$, that
	\beq
	\lim_{y\to0^{+}}\Omega_{2}(c+\ii y)m_{\nu}(\Omega_{2}(c+\ii y))=0,
	\eeq
	and hence
	\beq
	\lim_{y\to0^{+}}m_{\nu}(\Omega_{2}(c+\ii y))=0.
	\eeq
	On the other hand, as $\mu$ has a single interval support with density which is strictly positive in its interior, the absolute value of the Stieltjes transform $m_{\nu}$ is bounded below in $\caD$. Therefore we obtain a contradiction.
	
	Now suppose that $\Omega_{1}(c)=\infty$. As $\mu$ has mean $1$ with compact support, we have 
	\beq
	\lim_{y\to 0^{+}}\frac{\Omega_{2}(c+\ii y)}{c+\ii y}
	=\lim_{y\to 0^{+}}\frac{M_{\mu}(\Omega_{1}(c+\ii y))}{\Omega_{1}(c+\ii y)}
	=\lim_{z\to \infty}\frac{M_{\mu}(z)}{z}=1.
	\eeq
	Using Proposition~\ref{prop:newsubor} we have
	\beq
	\lim_{y\to0^{+}}M_{\nu}(\Omega_{2}(c+\ii y))=\lim_{y\to0^{+}}M_{\rho}(c+\ii y)=\lim_{y\to0^{+}}\frac{\Omega_{1}(c+\ii y)\Omega_{2}(c+\ii y)}{c+\ii y}=\infty,
	\eeq
	so that
	\beq
	\lim_{y\to0^{+}}\Omega_{2}(c+\ii y)m_{\nu}(\Omega_{2}(c+\ii y))+1=\lim_{y\to0^{+}}\left(1-M_{\nu}(\Omega_{2}(c+\ii y))\right)^{-1}=0.
	\eeq
	Since the measure $x\dd\mu(x)$ has strictly positive density in the interior of its support, $zm_{\nu}(z)+1$ is bounded below for $z\in\caD$, which is a contradiction.
\end{proof}

Recall the definitions of $\wh{I}_{\mu}$ and $\wh{I}_{\nu}$ in Definition \ref{def:I}. In the proof of Proposition \ref{prop:stab}, we need upper and lower bounds of the integrals $\wh{I}_{\mu}\circ\Omega_{1}$ and $\wh{I}_{\nu}\circ\Omega_{2}$. The following inequality together with Lemma \ref{lem:suborbd} will prove the desired bounds. 
\begin{lemma}\label{lem:bdeq}
	For any $z\in\C_{+}$, we have
	\beq\label{eq:hatIbdd}
	\absv{z}^{2}\wh{I}_{\mu}(\Omega_{1}(z))\wh{I}_{\nu}(\Omega_{2}(z))< 1.
	\eeq
\end{lemma}
\begin{proof}
	Let $z\in\C_{+}$ be fixed and we consider the following equality which follows from Remark \ref{rem:Iim}:
	\beq
	\absv{z}^{2}\wh{I}_{\mu}(\Omega_{1}(z))\wh{I}_{\nu}(\Omega_{2}(z))-1 =\frac{\absv{z}^{2}\im(\Omega_{1}(z)/z)\im(\Omega_{2}(z)/z)-\im \Omega_{1}(z)\im \Omega_{2}(z)}{\im \Omega_{1}(z)\im \Omega_{2}(z)}.
	\eeq
	Denoting $\arg z=\theta,\arg\Omega_{1}(z)=\theta_{1}$, and $\arg\Omega_{2}(z)=\theta_{2}$, the numerator is equal to the following after dividing by $\absv{\Omega_{1}(z)\Omega_{2}(z)}$:
	\begin{multline}
	\sin(\theta_{1}-\theta)\sin(\theta_{2}-\theta)-\sin\theta_{1}\sin\theta_{2}	=\frac{\cos(\theta_{1}+\theta_{2})-\cos(\theta_{1}+\theta_{2}-2\theta)}{2}	\\
	=\frac{\cos((\theta_{1}+\theta_{2}-\theta)+\theta)-\cos((\theta_{1}+\theta_{2}-\theta)-\theta)}{2}\\
	=\frac{1}{2}\left(\cos(\theta_{1}+\theta_{2}-\theta)\cos\theta-\sin(\theta_{1}+\theta_{2}-\theta)\sin\theta-\cos(\theta_{1}+\theta_{2}-\theta)\cos\theta -\sin(\theta_{1}+\theta_{2}-\theta)\sin\theta\right)	\\
	=-\sin\theta\sin(\theta_{1}+\theta_{2}-\theta).
	\end{multline}
	From Proposition \ref{prop:newsubor} we have
	\beq
	\theta_{1}+\theta_{2}-\theta=\arg(\Omega_{1}(z)\Omega_{2}(z)/z)=\arg(M_{\rho}(z)),
	\eeq
	so that
	\beq
	\absv{z}^{2}\wh{I}_{\mu}(\Omega_{1}(z))\wh{I}_{\nu}(\Omega_{2}(z))-1=-\frac{\absv{zM_{\rho}(z)}}{\im\Omega_{1}(z)\im\Omega_{2}(z)}\sin\theta\sin(\theta_{1}+\theta_{2}-\theta)=-\frac{\im z\im M_{\rho}(z)}{\im\Omega_{1}(z)\im\Omega_{2}(z)}<0,
	\eeq
	since $z,\Omega_{1}(z),\Omega_{2}(z)$, and $M_{\rho}(z)$ are all in $\C_{+}$.
\end{proof}

Combining Lemmas \ref{lem:suborbd} and \ref{lem:bdeq}, we can prove the following assertion:
\begin{lemma}\label{lem:intstab}
	Let $\mu$ and $\nu$ satisfy the assumptions of Theorem \ref{thm:main} and let $\caD\subset\ol{\C_{+}}\setminus\{0,\infty\}$ be compact. Then there exist constants $c_{1},c_{2},c_{3}>0$ such that the following hold:
	\begin{align}
	\inf_{z\in\caD}I_{\mu}(\Omega_{1}(z))\geq c_{1}, & & \inf_{z\in\caD}I_{\nu}(\Omega_{2}(z))\geq c_{1}, \label{eq:Ibd1}\\
	\inf_{z\in\caD}\wh{I}_{\mu}(\Omega_{1}(z))\geq c_{2}, & & \inf_{z\in\caD}\wh{I}_{\nu}(\Omega_{2}(z))\geq c_{2}, \label{eq:Ibd2}\\ 
	\sup_{z\in\caD}\wh{I}_{\mu}(\Omega_{1}(z))\leq c_{3}, & & \sup_{z\in\caD}\wh{I}_{\nu}(\Omega_{2}(z))\leq c_{3}. \label{eq:Ibd3}
	\end{align}
\end{lemma}
\begin{proof}
	The inequalities \eqref{eq:Ibd1} and \eqref{eq:Ibd2} follow directly from \eqref{eq:NPrepM} and  Lemma~\ref{lem:suborbd}. Given \eqref{eq:Ibd1} and \eqref{eq:Ibd2}, the inequality \eqref{eq:Ibd3} follows from Lemma \ref{lem:bdeq} by merely noting that $\absv{z}$ is bounded below and above in $\caD$. 
\end{proof}

The last estimate required to prove Proposition \ref{prop:stab} is the following purely computational lemma:
\begin{lemma}[Lemma 3.4 of \cite{Bao-Erdos-Schnelli2017arXiv}]\label{lem:comput}
	Let $z=E+\ii y$ with $y\geq 0$ and $\absv{z}\leq \theta$ for some small $\theta>0$. For $-1<t<1$, the following holds:
	\beq
	\int_{0}^{\theta}\frac{x^{t}}{(x-E)^{2}+y^{2}}\dd x\sim
	\begin{cases}
		\frac{E^{t}}{y}&\text{ if }E>y,\\
		\absv{z}^{t-1}\sim \absv{E}^{t-1}&\text{ if }  E<-y, \\
		y^{t-1} & \text{ if }y>\absv{E},
	\end{cases}
	\eeq
	where we write $C(z)\sim D(z)$ whenever there exists a constant $c>1$ such that $D(z)/c<C(z)<cD(z)$ uniformly in $\{z=E+\ii y:y>0,\absv{z}\leq\theta\}$.
\end{lemma}

Finally we can prove the main result of this subsection:
\begin{proposition}\label{prop:stab}
	Let $\mu$ and $\nu$ satisfy the assumptions of Theorem~\ref{thm:main}. Then there exists a constant $c>0$ such that
	\beq
	\inf_{z\in\C_{+}}\dist(\Omega_{1}(z),\supp\mu)\geq c,\quad 
	\inf_{z\in\C_{+}}\dist(\Omega_{2}(z),\supp\nu)\geq c.
	\eeq
\end{proposition}
\begin{proof}
	For simplicity, for $\omega\in\C_{+}$ we denote
	\beq
	d_{\mu}(\omega)= \dist(\omega,\supp \mu),\quad d_{\nu}(\omega)= \dist(\omega,\supp\nu).
	\eeq
	If $z$ is sufficiently close to $0$ or $\infty$, Lemma~\ref{lem:suborntlim} and Proposition~\ref{prop:suborext0} readily prove the result since $\mu$ and $\nu$ satisfy the assumptions of Theorem~\ref{thm:bound}. Thus, we may restrict our attention to a compact subset $\caD$ of $\ol{\C_{+}}\setminus\{0,\infty\}$ and it is enough to prove $\dist(\Omega_{1}(\caD),\supp\mu)>0$.
	
	By Lemmas~\ref{lem:suborext} (i) and \ref{lem:intstab}, it suffices to prove the following: if a sequence $\{\omega_{n}\}_{n\in\N}\in\Omega_{1}(\caD)\cap\C_{+}$ converges to a point $E\in\supp\mu$, then $\wh{I}_{\mu}(\omega_{n})$ diverges.  In order to do so, we first assume that $\omega_{n}$ converges to a point in $[E_{-}^{\mu},E_{-}^{\mu}+\delta]$ for some fixed $\delta>0$ to be chosen. 
	
	Recall that
	\begin{multline}\label{eq:hatIquot}
	\wh{I}_{\mu}(z)=\frac{\im\frac{M_{\mu}(z)}{z}}{\im z}
	=\frac{1}{\im z}\im\left[\frac{1}{z}-\frac{1}{z(zm_{\mu}(z)+1)}\right]
	=-\frac{1}{\absv{z}^{2}}+\frac{\im(zm_{\mu}(z)+1)}{(\im z)\absv{z}^{2}\absv{zm_{\mu}(z)+1}^{2}}	\\
	=\frac{1}{\absv{z}^{2}}\left[-1+\absv{\int_{\R_{+}}\frac{x}{x-z}\dd\mu(x)}^{-2}\left(\int_{\R_{+}}\frac{x}{\absv{x-z}^{2}}\dd\mu(x)\right)\right].
	\end{multline}
	Observe that if $z\in\C_{+}\cup\R\setminus\supp\mu$ satisfies $\absv{z-E^{\mu}_{-}}\leq\delta$ for sufficiently small $\delta$, we have
	\beq\label{eq:0thder}
	\absv{zm_{\mu}(z)+1}\leq C+C'\int_{0}^{2\delta}\frac{x^{t_{-}^{\mu}}}{\absv{x+E_{-}^{\mu}-z}}\dd x\leq C+C'(d_{\mu}(z))^{t_{-}^{\mu}}.
	\eeq
	Similarly, for $\absv{z-E_{-}^{\mu}}\leq\delta$, Lemma~\ref{lem:comput} gives us
	\beq\label{eq:1stder}
	\int_{\R_{+}}\frac{x\dd\mu(x)}{\absv{x-z}^{2}}\geq 
	\begin{cases}
		c\frac{(E-E_{-}^{\mu})^{t_{-}^{\mu}}}{y}&\text{ if }E-E_{-}^{\mu}>y,\\
		c(E_{-}^{\mu}-E)^{t_{-}^{\mu}-1}&\text{ if }  E-E_{-}^{\mu}<-y, \\
		cy^{t_{-}^{\mu}-1} & \text{ if }y>\absv{E-E_{-}^{\mu}},
	\end{cases}
	\eeq
	where we denote $z=E+\ii y$. 
	
	Then in each case of $t_{-}^{\mu}$ being nonnegative or negative, we analyze the quotient
	\beq\label{eq:quot}
	\left(\int_{\R_{+}}\frac{x}{\absv{x-z}^{2}}\dd\mu(x)\right)/\absv{\int_{\R_{+}}\frac{x}{x-z}\dd\mu(x)}^{2}.
	\eeq
	If $t_{-}^{\mu}\geq 0$, the RHS of \eqref{eq:0thder} is bounded, so that the quotient in \eqref{eq:quot} is bounded below by
	\beq
	\begin{cases}
		c\frac{(E-E_{-}^{\mu})^{t_{-}^{\mu}}}{y}&\text{ if }	E-E_{-}^{\mu}>y,\\
		c(E_{-}^{\mu}-E)^{t_{-}^{\mu}-1}			&\text{ if } 	 E-E_{-}^{\mu}<-y, \\
		cy^{t_{-}^{\mu}-1}							&\text{ if }	y>\absv{E-E_{-}^{\mu}},
	\end{cases}
	\eeq
	On the other hand if $t_{-}^{\mu}<0$, the last quantity in \eqref{eq:0thder} diverges with order $(d_{\mu}(z))^{t_{-}^{\mu}}$ as $z$ approaches to $\supp\mu$. In particular, from
	\beq
	d_{\mu}(z)\leq 
	\begin{cases}
		2(E-E_{-}^{\mu})		&\text{ if }E-E_{-}^{\mu}>y,\\
		2\absv{E-E_{-}^{\mu}}	&\text{ if }  E-E_{-}^{\mu}<-y, \\
		2y					&\text{ if }y>\absv{E-E_{-}^{\mu}},
	\end{cases}
	\eeq
	we conclude that the quotient in \eqref{eq:quot} has following lower bound:
	\beq
	\begin{cases}
		c\frac{(E-E_{-}^{\mu})^{-t_{-}^{\mu}}}{y}&\text{ if }	E-E_{-}^{\mu}>y,\\
		c(E_{-}^{\mu}-E)^{-t_{-}^{\mu}-1}			&\text{ if } 	 E-E_{-}^{\mu}<-y, \\
		cy^{-t_{-}^{\mu}-1}							&\text{ if }	y>\absv{E-E_{-}^{\mu}}.
	\end{cases}
	\eeq
	We see that in both cases, the quotient in \eqref{eq:quot} diverges as if we replace $z=\omega_{n}$. Since we have $\absv{\omega_{n}}\geq c$ from $0\notin\caD$ and Lemma \ref{lem:suborbd}, $\wh{I}_{\mu}(\omega_{n})$ also diverges by \eqref{eq:hatIquot}. Similar reasoning proves the same result when $\omega_{n}$ converges to a point in $[E_{+}^{\mu}-\delta,E_{+}^{\mu}]$. 
	
	Finally we consider the case in which $\omega_{n}$ converges to a point in the bulk of $\supp \mu$, namely $[E_{-}^{\mu}+\delta,E_{+}^{\mu}-\delta]$. Take $\delta'\in(0,\delta)$. Since the density $f_{\mu}$ of $\mu$ is strictly positive on $[E_{-}^{\mu}+\delta',E_{+}^{\mu}-\delta']$, there exists a constant $C>0$ depending on $\delta'$ such that
	\beq
	\im (z m_{\mu}(z)+1)=y\int_{\R_{+}}\frac{x\dd\mu(x)}{\absv{x-z}^{2}}\geq C,\quad \forall z\in\left(\C_{+}\cup\R\right)\setminus\supp\mu\cap\{z:\dist(z,[E_{-}^{\mu}+\delta,E_{+}^{\mu}-\delta])<\delta'\}.
	\eeq
	Also, for density being bounded in the bulk, there also exist constants $C_{1},C_{2}>0$ depending on $\delta'$ such that
	\beq
	\absv{z m_{\mu}(z)+1}
	\leq C_{1}+\int_{E_{-}^{\mu}+\delta}^{E_{+}^{\mu}-\delta}\frac{\dd x}{\absv{x-z}}
	\leq C_{1}+C_{2}\int_{E_{-}^{\mu}+\delta}^{E_{+}^{\mu}-\delta}\frac{\dd x}{\absv{x-E}+y}\leq C_{1}+C_{2}\absv{\log y}
	\eeq
	for all $z\in\C_{+}\cup\R\setminus\supp\mu$ with $\dist(z,[E_{-}^{\mu}+\delta,E_{+}^{\mu}-\delta])<\delta'$. Thus the quotient in \eqref{eq:quot} with $z=\omega_{n}$ diverges to infinity as $n\to\infty$, and so does $\wh{I}_{\mu}(\omega_{n})$. Therefore $\omega_{n}$ converging to $\supp\mu$ implies $\wh{I}_{\mu}(\omega_{n})\to\infty$, concluding the proof.
\end{proof}

\begin{remark}\label{rem:easecond}
	As easily seen, the power law behavior of $\mu$ and $\nu$ are used only in the proof of Proposition~\ref{prop:stab}. We remark that the same proof can be applied even if the power laws are replaced by the divergence of 
	\beq
	\absv{\int_{\R_{+}}\frac{x}{x-z}\dd\mu(x)}^{-2}\left(\int_{\R_{+}}\frac{x}{\absv{x-z}^{2}}\dd\mu(x)\right)
	\eeq
	for $z$ tending to the support of $\mu$ and the same for $\nu$.
\end{remark}

\begin{remark}\label{rem:0stab}
	If $E_{-}^{\mu}$ or $E_{-}^{\nu}$ is $0$, we can still prove that $\Omega_{1}(z)$ and $\Omega_{2}(z)$ stay way from the support if $z$ is away from $0$. On the other hand, both of the subordination functions must tend to $0$ as $z$ tends to $0$ nontangentially(see Lemma \ref{lem:suborntlim}).
\end{remark}

\subsection{Proof of Theorem \ref{thm:main}}
Given Proposition \ref{prop:stab} we proceed to the next step, to prove that the density $f$ is strictly positive on a single interval and $0$ elsewhere(see Proposition \ref{prop:edgechar2} for precise statements). The first thing to do is to observe that $\im\Omega_{1}(E)$, $\im\Omega_{2}(E)$, and $f(E)$ are all comparable for $E\in\R$. The density $f(E)$ is comparable to $\im\Omega_{1}(E)$ using
\beq\label{eq:dens_imsubor_compa}
Ef(E)=\lim_{y\to 0^{+}}\im(zm_{\rho}(z)+1)=\lim_{y\to 0^{+}}\im(\Omega_{1}(z)m_{\mu}(z)+1)=\im\Omega_{1}(E) I_{\mu}(\Omega_{1}(E)),
\eeq
where we denoted $z=E+\ii y$. We also have the same equality with $(\mu,\nu)$ and $(\Omega_{1},\Omega_{2})$ interchanged, so that
\beq\label{eq:imsubor_compa}
\lim_{y\to0^{+}}\frac{\im\Omega_{2}(z)}{\im\Omega_{1}(z)}=\frac{I_{\mu}(\Omega_{1}(E))}{I_{\nu}(\Omega_{2}(E))}.
\eeq
Note that by Lemma \ref{lem:suborext} and Proposition \ref{prop:stab}, we can take the limits in \eqref{eq:dens_imsubor_compa} and \eqref{eq:imsubor_compa} in $\C_{+}$ (not only vertically), and the functions $E\mapsto I_{\mu}(\Omega_{1}(E)),I_{\nu}(\Omega_{2}(E))$ are continuous and strictly positive.

The proof of Proposition \ref{prop:edgechar2} is again divided into two steps: we first characterize the \emph{edges} of $\supp(\mu\boxplus\nu)$, which is defined as follows:
\beq
\caV= \partial\{x\in\R:f(x)>0\}.
\eeq
After characterizing $\caV$, we prove that there are exactly two points that match the characterization to conclude. The following result corresponds to the first step, characterization of the edges.
\begin{proposition}\label{prop:edgechar}
	Let $\mu$ and $\nu$ satisfy the assumptions of Theorem \ref{thm:main}. Then for all $z\in\C_{+}\cup\R$, the following inequality holds:
	\beq\label{eq:edgecharineq}
	\absv{\left(\frac{\Omega_{1}(z)}{M_{\mu}(\Omega_{1}(z))}M_{\mu}'(\Omega_{1}(z))-1\right)\left(\frac{\Omega_{2}(z)}{M_{\nu}(\Omega_{2}(z))}M_{\nu}'(\Omega_{2}(z))-1\right)}\leq1.
	\eeq
	Furthermore, for $z=E+\ii y\in\C_{+}\cup\R$, the equality holds in \eqref{eq:edgecharineq} if and only if $E\in\caV$ and $y=0$. In this case, we also have
	\beq\label{eq:edgechar}
	\left(\frac{\Omega_{1}(z)}{M_{\mu}(\Omega_{1}(z))}M_{\mu}'(\Omega_{1}(z))-1\right)
	\left(\frac{\Omega_{2}(z)}{M_{\nu}(\Omega_{2}(z))}M_{\nu}'(\Omega_{2}(z))-1\right)=1.
	\eeq
\end{proposition}

Before the proof, we recall that the images of $\Omega_{1}$ and $\Omega_{2}$ are away from $\supp\mu$ and $\supp \nu$, respectively. As a result, the compositions $M_{\mu}^{(\ell)}\circ \Omega_{1}$ and $M_{\nu}^{(\ell)}\circ\Omega_{2}$ are well-defined analytic functions on $\C_{+}$ which extend continuously to $\R$ by Lemma \ref{lem:suborext} and Proposition \ref{prop:suborext0}, where we denote $M_{\mu}^{(\ell)}=\frac{\dd^{\ell}}{\dd z^{\ell}}M_{\mu}$. The same applies to $H_{\mu}\circ\Omega_{1}$ and $H_{\nu}\circ\Omega_{2}$ by Lemma \ref{lem:NPrepMsupp}. Furthermore, the denominators $M_{\mu}\circ\Omega_{1}$ and $M_{\nu}\circ\Omega_{2}$ in \eqref{eq:edgecharineq} and \eqref{eq:edgechar} are nonzero on $\C_{+}\cup\R$ by Proposition \ref{prop:newsubor} (iv) and Lemma \ref{lem:suborbd}.

\begin{proof}[Proof of Proposition \ref{prop:edgechar}]
	We first observe from Lemma~\ref{lem:NPrepM} that
	\beq
	M_{\mu}'(z)=\frac{\dd}{\dd z}(zm_{\wh{\mu}}(z)+z)=\int_{\R_{+}}\frac{x}{(x-z)^{2}}\dd\wh{\mu}+1,
	\eeq
	so that
	\begin{multline}
	\frac{\Omega_{1}(z)}{M_{\mu}(\Omega_{1}(z))}M_{\mu}'(\Omega_{1}(z))-1
	=\frac{1}{M_{\mu}(\Omega_{1}(z))}\int_{\R_{+}}\left(\frac{x\Omega_{1}(z)}{(x-\Omega_{1}(z))^{2}}\dd\wh{\mu}(x)+\Omega_{1}(z)-M_{\mu}(\Omega_{1}(z))\right)\dd\wh{\mu}(x)\\
	=\frac{1}{M_{\mu}(\Omega_{1}(z))}\int_{\R_{+}}\left[\frac{x\Omega_{1}(z)}{(x-\Omega_{1}(z))^{2}}-\frac{\Omega_{1}(z)}{x-\Omega_{1}(z)}\right]\dd\wh{\mu}(x)
	=\frac{\Omega_{1}(z)^{2}}{M_{\mu}(\Omega_{1}(z))}\int_{\R_{+}}\frac{1}{(x-\Omega_{1}(z))^{2}}\dd\wh{\mu}(x).
	\end{multline}
	
	Thus, by Lemma \ref{lem:bdeq} we have
	\begin{multline}\label{eq:ineqchain}
	\absv{\left(\frac{\Omega_{1}(z)}{M_{\mu}(\Omega_{1}(z))}M_{\mu}'(\Omega_{1}(z))-1\right)\left(\frac{\Omega_{2}(z)}{M_{\nu}(\Omega_{2}(z))}M_{\nu}'(\Omega_{2}(z))-1\right)}	\\
	\overset{(i)}{\leq}\absv{\frac{\Omega_{1}(z)^{2}\Omega_{2}(z)^{2}}{M_{\rho}(z)^{2}}}\wh{I}_{\mu}(\Omega_{1}(z))\wh{I}_{\nu}(\Omega_{2}(z))
	=\absv{z}^{2}\wh{I}_{\mu}(\Omega_{1}(z))\wh{I}_{\nu}(\Omega_{2}(z))\overset{(ii)}{\leq} 1,
	\end{multline}
	concluding the proof for the inequality \eqref{eq:edgecharineq}.
	
	Now we prove that equality in \eqref{eq:edgecharineq} implies $z=E\in\caV$. Let $z=E+\ii y\in\C_{+}$ be such that the equality holds in \eqref{eq:edgecharineq} and hence inequalities (\textit{i}) and (\textit{ii}) in \eqref{eq:ineqchain} are both equalities. First, (\textit{i}) being an equality implies
	\beq
	\absv{\int_{\R_{+}}\frac{1}{(x-\Omega_{1}(z))^{2}}\dd\wh{\mu}(x)}=\int_{\R_{+}}\frac{1}{\absv{x-\Omega_{1}(z)}^{2}}\dd\wh{\mu}(x),
	\eeq
	so that $(x-\Omega_{1}(z))^{2}=\absv{x-\Omega_{1}(z)}^{2}$ for $\wh{\mu}-$a.e. $x\in\R$. In particular, for $\absv{\Omega_{1}(z)}$ being finite by Lemma \ref{lem:suborbd}, we have $\im\Omega_{1}(z)=0$. By a similar reasoning we also have $\im\Omega_{2}(z)=0$.
	
	Then, noting that equality in (\textit{ii}) prevents $z=0$, we have
	\beq
	\frac{\Omega_{2}(z)}{z}=\frac{M_{\mu}(\Omega_{1}(z))}{\Omega_{1}(z)}
	=\frac{\int \frac{1}{x-\Omega_{1}(z)}\dd\mu(x)}{\int_{\R_{+}}\frac{x}{x-\Omega_{1}(z)}\dd\mu(x)}>0,
	\eeq
	and similarly
	\beq
	\frac{\Omega_{1}(z)}{z}=\frac{M_{\nu}(\Omega_{2}(z))}{\Omega_{2}(z)}>0,
	\eeq
	using Proposition~\ref{prop:stab} and $\im\Omega_{2}(z)=\im \Omega_{1}(z)=0$. Thus $z,\Omega_{1}(z),\Omega_{2}(z)$, and $M_{\rho}(z)$ are nonzero real numbers with the same sign.
	
	If $z<0$ so that $\Omega_{1}(z),\Omega_{2}(z)$, and $M_{\rho}(z)$ are all negative, from the definition of $M_{\mu}$ we find that
	\begin{multline}
	\frac{\Omega_{1}(z)}{M_{\mu}(\Omega_{1}(z))}M_{\mu}'(\Omega_{1}(z))	=\frac{1}{(\Omega_{1}(z)m_{\mu}(\Omega_{1}(z))+1)m_{\mu}(\Omega_{1}(z))}\int_{\R_{+}}\frac{x}{(x-\Omega_{1}(z))^{2}}\dd\mu(x)	\\
	=\left(\int_{\R_{+}}\frac{x}{x-\Omega_{1}(z)}\dd\mu(x)\right)^{-1}\left(\int_{\R_{+}}\frac{1}{x-\Omega_{1}(z)}\dd\mu(x)\right)^{-1}\left(\int_{\R_{+}}\frac{x}{(x-\Omega_{1}(z))^{2}}\dd\mu(x)\right)\in (0,1],
	\end{multline}
	where the upper bound follows from Cauchy-Schwarz inequality applied to two functions $(x/(x-z))^{1/2}$ and $(1/(x-z))^{1/2}$. Thus both of the factors in \eqref{eq:edgecharineq} are strictly less than $1$ in modulus, so that the equality is not possible for $z<0$.
	
	Given the fact that $z=E,\Omega_{1}(E)$, and $\Omega_{2}(E)$ are all positive, we see from \eqref{eq:imsubor_compa} that
	\begin{multline}\label{eq:edgechar_pf1}
	0<E\wh{I}_{\mu}(\Omega_{1}(E))
	=\lim_{\theta\to 0^{+}}E\frac{\im\frac{\Omega_{2}(E\e{\ii\theta})}{E\e{\ii\theta}}}{\im \Omega_{1}(E\e{\ii\theta})}	
	=\lim_{\theta\to 0^{+}}\frac{\im\Omega_{2}(E\e{\ii\theta})}{\im \Omega_{1}(E\e{\ii\theta})}\frac{E\im\frac{\Omega_{2}(E\e{\ii\theta})}{E\e{\ii\theta}}}{\im\Omega_{2}(E\e{\ii\theta})}
	=\frac{I_{\mu}(\Omega_{1}(E))}{I_{\nu}(\Omega_{2}(E))}\lim_{\theta\to 0^{+}}\frac{\sin(\theta_{2}-\theta)}{\sin\theta_{2}}
	\end{multline}
	and similarly that
	\beq\label{eq:edgechar_pf2}
	E\wh{I}_{\nu}(\Omega_{2}(E))=\frac{I_{\nu}(\Omega_{2}(E))}{I_{\mu}(\Omega_{1}(E))}\lim_{\theta\to 0^{+}}\frac{\sin(\theta_{1}-\theta)}{\sin\theta_{1}},
	\eeq
	where we denoted $\arg\Omega_{1}(E\e{\ii\theta})=\theta_{1}$ and $\arg\Omega_{2}(E\e{\ii\theta})=\theta_{2}$. Note that the quantities $E\wh{I}_{\mu}(\Omega_{1}(E))$, $E\wh{I}_{\nu}(\Omega_{2}(E))$, and ${I_{\mu}(\Omega_{1}(E))}/{I_{\nu}(\Omega_{2}(E))}$ are all positive, so that
	\begin{align}
	&\lim_{\theta\to 0^{+}}\frac{\sin(\theta_{1}-\theta)}{\sin\theta_{1}}>0,
	&&\lim_{\theta\to 0^{+}}\frac{\sin(\theta_{2}-\theta)}{\sin\theta_{2}}>0.
	\end{align}
	Also $\sin\theta_{1}$ and $\sin\theta_{2}$ must converge to $0$, since $\Omega_{1}(z)$ and $\Omega_{2}(z)$ converge to positive numbers $\Omega_{1}(E)$ and $\Omega_{2}(E)$ as $\theta$ tends to $0$. Multiplying the equalities yields
	\beq
	1=E^{2}\wh{I}_{\mu}(\Omega_{1}(E))\wh{I}_{\nu}(\Omega_{2}(E))=\left(\lim_{\theta\to 0^{+}}\frac{\sin(\theta_{1}-\theta)}{\sin\theta_{1}}\right)\left(\lim_{\theta\to 0^{+}}\frac{\sin(\theta_{2}-\theta)}{\sin\theta_{2}}\right),
	\eeq
	where the first equality is exactly (\textit{ii}). Since $\theta_{1},\theta_{2}\geq\theta$ by Proposition \ref{prop:newsubor} (ii) and $\theta_{1},\theta_{2}\to0$, we get
	\beq
	\lim_{\theta\to 0^{+}}\frac{\sin(\theta_{1}-\theta)}{\sin\theta_{1}}
	=\lim_{\theta\to 0^{+}}\frac{\theta_{1}-\theta}{\theta_{1}}=1
	=\lim_{\theta\to 0^{+}}\frac{\sin(\theta_{2}-\theta)}{\sin\theta_{2}}
	=\lim_{\theta\to 0^{+}}\frac{\theta_{2}-\theta}{\theta_{2}}.
	\eeq
	To summarize, so far we have proved that for any $z=E+\ii y\in\C_{+}\cup\R$ at which the equality in \eqref{eq:edgecharineq} holds, all three of $z=E$, $\Omega_{1}(E),\Omega_{2}(E)$ must be positive real numbers and
	\begin{align}
	\lim_{\theta\to 0^{+}}\arg\Omega_{1}(E\e{\ii\theta})&=0=\lim_{\theta\to 0^{+}}\arg\Omega_{2}(E\e{\ii\theta}),&
	\lim_{\theta\to 0^{+}}\frac{\theta}{\arg \Omega_{1}(E\e{\ii\theta})}&=0=\lim_{\theta\to 0^{+}}\frac{\theta}{\arg \Omega_{2}(E\e{\ii\theta})}.
	\end{align}
	
	We see that $\im \Omega_{2}(E)=0$ and $I_{\nu}(\Omega_{2}(E))<\infty$ implies $\im (Em_{\rho}(E)+1)=0$, and thus $Ef(E)=0$. Now if we suppose that $E\in (\ol{\{x\in\R:f(x)>0\}})^{c}$, or equivalently, $\dist(E,\supp\rho)>0$, then for a small enough $\theta_{0}>0$ there exists a constant $C>0$ such that
	\beq
	\frac{\im (E\e{\ii\theta}m_{\rho}(E\e{\ii\theta})+1)}{E\sin\theta}=\frac{\im (E\e{\ii\theta}m_{\rho}(E\e{\ii\theta})+1)}{\im E\e{\ii\theta}}=\int_{\R_{+}}\frac{x}{\absv{x-E\e{\ii\theta}}^{2}}\dd\rho(x)\leq C,
	\eeq
	for all $\theta\in[0,\theta_{0}]$. This in turn implies 
	\beq
	\frac{\sin \theta}{\sin\theta_{2}}=\frac{\absv{\Omega_{2}(E\e{\ii\theta})}}{E}\frac{E\sin\theta}{\im (E\e{\ii\theta}m_{\rho}(E\e{\ii\theta})+1)}\frac{\im (E\e{\ii\theta}m_{\rho}(E\e{\ii\theta})+1)}{\im \Omega_{2}(E\e{\ii\theta})}	
	\geq \frac{1}{C}\frac{\absv{\Omega_{2}(E\e{\ii\theta})}}{E}I_{\nu}(\Omega_{2}(E))\geq c,
	\eeq
	for all $\theta\in[0,\theta_{0}]$, contradicting $\theta/\theta_{2}\to 0$. 
	
	For the converse, we take $z=E\in \caV$ and will prove that (\textit{i}) and (\textit{ii}) in \eqref{eq:ineqchain} are equalities. As the density is continuous and bounded by Theorem~\ref{thm:bound}, we readily see that $f(E)=0$ implies $\im m_{\rho}(E)=0$ and thus $\im (Em_{\rho}(E)+1)=0$. Furthermore, since $\rho$ is compactly supported in $(0,\infty)$, $\caV$ must also be contained in $(0,\infty)$ so that $E>0$. Recalling
	\beq
	\im \Omega_{2}(E)I_{\nu}(\Omega_{2}(E))=\im(Em_{\rho}(E)+1)
	\eeq
	we also get $\im \Omega_{2}(E)=\im \Omega_{1}(E)=0$, proving that (\textit{i}) is an equality. As before $H_{\mu}(\Omega_{1}(E))=\Omega_{2}(E)/E>0$ implies that $\Omega_{1}(E),\Omega_{2}(E)$ have the same sign as $E$, and hence are positive.
	
	Now we prove that equality holds in (\textit{ii}). By the definition of $\caV$, there must be a sequence $\{\epsilon_{n}\}_{n\in\N}$ of real numbers increasing or decreasing to $0$ such that $\rho(E+\epsilon_{n})>0$ for all $n\in\N$. For each fixed $n\in\N$, Using the Stieltjes inversion and the continuity of $m_{\rho}$, we first observe that
	\begin{multline}
	0<(E+\epsilon_{n})\rho(E+\epsilon_{n})
	=\lim_{y\to 0^{+}}\im\left[(E+\epsilon_{n}+\ii y)m_{\rho}(E+\epsilon_{n}+\ii y)+1\right]\\
	=\lim_{\theta\to 0}\im \left[(E+\epsilon_{n})\e{\ii\theta}m_{\rho}((E+\epsilon_{n})\e{\ii\theta})+1\right]
	=\lim_{\theta\to 0}(E+\epsilon_{n})\sin\theta \int_{\R_{+}}\frac{x}{\absv{x-(E+\epsilon_{n})\e{\ii\theta}}^{2}}\dd \rho(x).
	\end{multline}
	Then Proposition~\ref{prop:stab} implies 
	\beq
	\lim_{\theta\to 0}\frac{\im \Omega_{2}((E+\epsilon_{n})
		\e{\ii\theta})}{(E+\epsilon_{n})\sin\theta}
	=\lim_{\theta\to0}\frac{1}{I_{\nu}(\Omega_{2}((E+\epsilon_{n})\e{\ii\theta}))}\int_{\R_{+}}\frac{x}{\absv{x-(E+\epsilon_{n})\e{\ii\theta}}^{2}}\dd\rho(x)=\infty,
	\eeq
	so that by denoting $\theta_{2}^{(n)}=\arg\Omega_{2}((E+\epsilon_{n})\e{\ii\theta})$ we get
	\beq
	0\leq\lim_{\theta\to 0}\frac{\theta}{\theta_{2}^{(n)}}
	\leq\lim_{\theta\to 0}\frac{\theta}{\sin\theta}\frac{\sin\theta}{\sin \arg \Omega_{2}((E+\epsilon_{n})\e{\ii\theta})}
	=\lim_{\theta\to 0} \frac{\sin \theta }{\sin \arg \Omega_{2}((E+\epsilon_{n})\e{\ii\theta})}=0,
	\eeq
	where we used $\absv{\Omega_{2}(E+\epsilon_{n})}<\infty$ in the last equality. Finally we conclude
	\beq
	\lim_{\theta\to 0}\frac{\sin(\theta_{2}^{(n)}-\theta)}{\sin \theta_{2}^{(n)}}=1-\lim_{\theta\to 0}\frac{\theta}{\theta_{2}^{(n)}}=1
	\eeq
	and the same equality for $\theta_{1}$, from which we deduce (\textit{ii}) using \eqref{eq:edgechar_pf1} and \eqref{eq:edgechar_pf2} and taking the limit $n\to\infty$.
	
	For the last part, we observe that $z$, $\Omega_{2}(z)$, and $\Omega_{1}(z)$ are all positive real numbers when the equivalent conditions hold, so that (\textit{i}) and (\textit{ii}) are equalities without absolute values on the left sides.
\end{proof}

Now that we are given Proposition~\ref{prop:edgechar}, we can characterize points in $\caV$ as the solutions of equation \eqref{eq:edgechar}. Yet we still do not know whether the set $\caV$ consists of exactly two points, $E_{-}$ and $E_{+}$. In order to exhaust the possibility of a non-edge point in $\caV$, that is, the existence of \emph{isolated zero} (see \cite{Bao-Erdos-Schnelli2018}), we need representations of the subordination functions corresponding to those of $H$-functions in Lemma~\ref{lem:NPrepM}.

\begin{lemma}\label{lem:NPrepSubor}
	Let $\mu$ and $\nu$ satisfy the assumptions of Theorem \ref{thm:main} and $\Omega_{1}$ and $\Omega_{2}$ be the subordination functions corresponding to $\mu$ and $\nu$ by means of Proposition~\ref{prop:newsubor}. Then there exists finite Borel measures $\wt{\mu}$ and $\wt{\nu}$ on $\R_{+}$ such that the following hold:
	\begin{itemize}
		\item[(i)] 
		$\frac{\Omega_{1}(z)}{z}=1+m_{\wt{\nu}}(z)$ and $\frac{\Omega_{2}(z)}{z}=1+m_{\wt{\mu}}(z)$ whenever $z\notin \supp\wt{\nu}$ and $z\notin\supp\wt{\mu}$, respectively,
		\item[(ii)] 
		$\wt{\mu}(\R)=\Var{\mu}$ and $\wt{\nu}(\R)=\Var{\nu}$,
		\item[(iii)] 
		$\supp\wt{\mu}=\supp\wt{\nu}=\supp \rho$.
	\end{itemize}
\end{lemma}
\begin{proof}
	We prove the lemma only for $\Omega_{2}$, and the proof for $\Omega_{1}$ is completely analogous. We start from the following identity:
	\beq
	\frac{\Omega_{2}(z)}{z}=H_{\mu}(\Omega_{1}(z)).
	\eeq
	Since $H_{\mu}$ maps $\C_{+}$ into itself as $\mu$ is nondegenerate, we see that $z\mapsto\Omega_{2}(z)/z$ is an analytic self-map of $\C_{+}$. Furthermore, we see from Lemma~\ref{lem:suborext} and \ref{lem:suborntlim} that $\Omega_{1}(\ii y)$ and $\Omega_{2}(\ii y)$ tend to infinity in $\C_{+}$ as $y\to\infty$ so that
	\beq
	\lim_{y\to\infty}\frac{\Omega_{2}(\ii y)}{\ii y}=\lim_{y\to\infty}H_{\mu}(\Omega_{1}(\ii y))=\lim_{z\to\infty,z\in\C_{+}}H_{\mu}(z)=1,
	\eeq
	where we used that the support of $\mu$ is compact. Similarly $\Omega_{1}(\ii y)/(\ii y)\to 1$ as $y\to\infty$.
	
	As $z\mapsto\Omega_{2}(z)/z$ is a Nevanlinna function on $\C_{+}$, similarly to Lemma \ref{lem:NPrepM} it suffices to prove that
	\beq
	\sup_{y\geq 1}\absv{\Omega_{2}(\ii y)-\ii y}<\infty.
	\eeq
	Indeed, by Lemma~\ref{lem:NPrepM} we see that 
	\beq
	\lim_{y\to\infty}\Omega_{2}(\ii y)-\ii y=\lim_{y\to\infty}\frac{\ii y}{\Omega_{1}(\ii y)}\left(M_{\mu}(\Omega_{1}(\ii y))-\Omega_{1}(\ii y)\right)=-\Var{\mu},
	\eeq
	which proves the first two assertions as $\Omega_{2}$ is continuous on the positive imaginary axis $\ii(0,\infty)$.
	
	In order to prove the last part, we first prove the inclusion $\supp\rho\subset\supp\wt{\mu}$. Let $U\subset\R\setminus\supp\wt{\mu}$ be an open interval, so that $\Omega_{2}$ is analytic and real-valued in $U$. Now it directly follows from \eqref{eq:dens_imsubor_compa} that $f$ is identically zero in $U$. Therefore we have the inclusion $\supp\rho\subset\supp\wt{\mu}$.
	
	For the converse, let $I\subset\R\setminus \supp\rho$ be an open interval. Then $f$ is identically zero in $I$, and so are $\im\Omega_{1}$ and $\im\Omega_{2}$ by \eqref{eq:dens_imsubor_compa}. Furthermore, for $E=0$ we have
	\beq
	\lim_{z\to 0}\frac{\Omega_{2}(z)}{z}=\lim_{z\to 0}H_{\mu}(\Omega_{1}(z))=\lim_{w\to0}H_{\mu}(w)=\int_{\R_{+}}\frac{1}{x}\dd\mu(x)\in\R.
	\eeq
	Thus (the analytic continuation of) $z\mapsto \Omega_{2}(z)/z$ has zero imaginary part in $I$. Therefore $I\subset\R\setminus\supp\wt{\mu}$, concluding the proof of $\supp \wt{\mu}\subset\supp\rho$. 
\end{proof}

Now we are ready to prove that $\caV$ in fact is exactly two endpoints of $\supp\rho$.
\begin{proposition}\label{prop:edgechar2}
	Let $\mu$ and $\nu$ satisfy the assumptions of Theorem \ref{thm:main}. Then there exist positive real numbers $E_{-}<E_{+}$ such that $\caV=\{E_{-},E_{+}\}$ and 
	\beq
	\{x\in\R:f(x)>0\}=(E_{-},E_{+}).
	\eeq
\end{proposition}
\begin{proof}
	We have seen that $\im\Omega_{1}(E)=0=\im \Omega_{2}(E)$ for any $E\in\caV$ in the proof of Lemma~\ref{prop:edgechar}. Thus using Proposition~\ref{prop:stab}, we divide the possible locations of $\Omega_{1}(E)$ and $\Omega_{2}(E)$ for $E\in\caV$ as follows:
	\begin{itemize}
		\item[(i)] 
		$\Omega_{1}(E)<E_{-}^{\mu}$ and $\Omega_{2}(E)<E_{-}^{\nu}$,
		\item[(ii)] 
		$\Omega_{1}(E)>E_{+}^{\mu}$ and $\Omega_{2}(E)>E_{+}^{\nu}$,
		\item[(iii)] 
		$\Omega_{1}(E)<E_{-}^{\mu}$ and $\Omega_{2}(E)>E_{+}^{\nu}$, or $\Omega_{1}(E)>E_{+}^{\mu}$ and $\Omega_{2}(E)<E_{-}^{\nu}$.
	\end{itemize}
	We will prove that in each case of (i) and (ii) the equation \eqref{eq:edgechar} have exactly one solution, while in the last case it does not have any. 
	
	For simplicity, recalling $\eqref{eq:edgechar}$, we define
	\beq
	h(E)= E^{2}\wh{I}_{\mu}(\Omega_{1}(E))\wh{I}_{\nu}(\Omega_{2}(E)), \quad \forall E>0.
	\eeq
	Note that the function $h$ has another form:
	\beq\label{eq:edgechar2}
	h(E)=\frac{1}{M_{\mu}(\Omega_{1}(E))M_{\nu}(\Omega_{2}(E))}\left(\int_{\R_{+}}\frac{\Omega_{1}(E)^{2}}{(x-\Omega_{1}(E))^{2}}\dd\wh{\mu}(x)\right)\left(\int_{\R_{+}}\frac{\Omega_{2}(E)^{2}}{(x-\Omega_{2}(E))^{2}}\dd\wh{\nu}(x)\right).
	\eeq
	Also note that Lemma~\ref{lem:NPrepSubor} implies
	\beq
	\Omega_{2}'(z)=\frac{\dd }{\dd z}(z+zm_{\wt{\mu}(z)})= 1+\int_{\R_{+}}\frac{x}{(x-z)^{2}}\dd\wt{\mu}(x),
	\eeq
	so that $\Omega_{2}$ is increasing on $(\supp\wt{\mu})^{c}=(\supp\rho) ^{c}$. Similarly $\Omega_{1}$ is increasing on $(\supp\rho)^{c}$. 
	
	We first prove the existence and uniqueness of $E$ in the cases (i) and (ii). For existence, let $E_{-}$ and $E_{+}$ be the leftmost and rightmost endpoints of $\supp\rho$, respectively. The existence of $E_{-}$ and $E_{+}$ follows from the fact that $f$ is a nonzero, continuous function on $\R$ with compact support. Clearly the points $E_{-}$ and $E_{+}$ are in $\caV$, so that they must solve \eqref{eq:edgechar} by Proposition \ref{prop:edgechar}. 
	
	By Lemma~\ref{lem:NPrepSubor}, we see that $\Omega_{2}$ and $\Omega_{1}$ map $(-\infty,E_{-})$ and $(E_{+},\infty)$ into real half lines with the same directions. Furthermore, by Proposition~\ref{prop:stab}, these images are exactly the leftmost and rightmost components of $\Omega_{2}(\R)\cap\R$ and $\Omega_{1}(\R)\cap\R$ so that \beq
	\Omega_{2}(E_{-})<E_{-}^{\nu}, \quad \Omega_{1}(E_{-})<E_{-}^{\mu}, \AND \Omega_{2}(E_{+})>E_{+}^{\nu}, \quad \Omega_{1}(E_{+})>E_{+}^{\mu}.
	\eeq
	
	To prove the uniqueness, suppose that $E_{0}$ satisfies (i) and \eqref{eq:edgechar}. If $E_{0}<E_{-}$, as $\wh{I}_{\mu}$ and $\wh{I}_{\nu}$ increase on $(-\infty,E_{-}^{\mu})$ and $(-\infty,E_{-}^{\nu})$ respectively, we see from $\Omega_{1}(E_{0})<\Omega_{1}(E_{-})$ and $\Omega_{2}(E_{0})<\Omega_{2}(E_{-})$ that
	\beq
	h(E_{0})=E_{0}^{2}\wh{I}_{\mu}(\Omega_{1}(E_{0}))\wh{I}_{\nu}(E_{0})<E_{-}^{2}\wh{I}_{\mu}(\Omega_{1}(E_{-}))\wh{I}_{\nu}(\Omega_{2}(E_{-}))=h(E_{-})=1,
	\eeq
	which is contradiction. On the other hand if $E_{0}$ is larger than $E_{-}$, we must have 
	\beq
	\wh{I}_{\mu}(\Omega_{1}(E_{0}))\wh{I}_{\nu}(\Omega_{2}(E_{0}))=\frac{1}{E_{0}^{2}}<\frac{1}{E_{-}^{2}}=\wh{I}_{\mu}(\Omega_{1}(E_{-}))\wh{I}_{\nu}(\Omega_{2}(E_{-})),
	\eeq
	so that either $\Omega_{1}(E_{0})<\Omega_{1}(E_{-})$ or $\Omega_{2}(E_{0})<\Omega_{2}(E_{-})$ must hold since $\wh{I}_{\mu}$ and $\wh{I}_{\nu}$ are increasing respectively on $(-\infty,E_{-}^{\mu})$ and $(-\infty,E_{-}^{\nu})$. Supposing the latter without loss of generality, it follows that $\Omega_{2}(E_{0})$ coincides with $\Omega_{2} (E_{1})$ for some $E_{1}\in(0,E_{-})$, as $\Omega_{2}$ is a continuous, strictly increasing function which maps $(-\infty,E_{-})$ onto $(-\infty,\Omega_{2}(E_{-}))$. Then we see that $\Omega_{1}(E_{0})$ must also be the same as $\Omega_{1}(E_{1})$, as both of them are the unique solution of $M_{\mu}(\Omega)=M_{\nu}(\Omega_{2}(E_{0}))=M_{\nu}(\Omega_{2}(E_{1}))$ in $(-\infty, E_{-}^{\mu})$. Now we observe from \eqref{eq:edgechar2} that $h(E_{0})=h(E_{1})<1$, which is a contradiction.
	
	Similarly, let $E_{0}\in\caV$ satisfy (ii). We first recall that $M_{\mu}$ and $M_{\nu}$ are positive and increasing on $(E_{+}^{\mu},\infty)$ and $(E_{+}^{\nu},\infty)$, and that
	\beq
	\frac{\dd}{\dd w}\frac{w^{2}}{(x-w)^{2}}=2\frac{w}{x-w}\frac{x}{(x-w)^{2}}<0,\quad \forall w>x>0.
	\eeq
	Therefore from \eqref{eq:edgechar2} we see that $h$ is decreasing on $(E_{+},\infty)$, so that $E_{0}>E_{+}$ implies $h(E_{0})<h(E_{+})=1$. On the other hand if $E_{0}<E_{+}$, we obtain 
	\beq
	\wh{I}_{\mu}(\Omega_{1}(E_{0}))\wh{I}_{\nu}(\Omega_{2}(E_{0})>\wh{I}_{\mu}(\Omega_{1}(E_{+}))\wh{I}_{\nu}(\Omega_{2}(E_{+})).
	\eeq
	Combining this inequality with (ii), either $\Omega_{2}(E_{0})<\Omega_{2}(E_{+})$ or $\Omega_{1}(E_{0})<\Omega_{1}(E_{+})$ must hold, and by Proposition~\ref{prop:newsubor} (iii) one implies the other using the same reasoning as in the case (i). Therefore we have
	\begin{align}
	E_{+}^{\mu}&<\Omega_{1}(E_{0})<\Omega_{1}(E_{+}),&
	E_{+}^{\nu}&<\Omega_{2}(E_{0})<\Omega_{2}(E_{+}).
	\end{align}
	By \eqref{eq:edgechar2} this would imply $h(E_{0})>h(E_{+})=1$, which is a contradiction,
	
	It remains to prove that there is no solution to \eqref{eq:edgechar} satisfying (iii). To this end, we suppose that $E_{0}$ is a solution satisfying
	\begin{align}\label{eq:edge1}
	\Omega_{1}(E_{0})&>E_{+}^{\mu},&
	\Omega_{2}(E_{0})&<E_{-}^{\nu}.
	\end{align}
	We note that \eqref{eq:edge1} readily implies $E_{0}\in(E_{-},E_{+})$. If not, say $E_{0}<E_{-}$, we would end up
	\begin{align}
	\Omega_{1}(E_{0})&\in(-\infty,\Omega_{1}(E_{-}))\subset(-\infty, E_{-}^{\mu}), &\Omega_{2}(E_{0})&\in(-\infty,\Omega_{2}(E_{-}))\subset(-\infty,E_{-}^{\nu}),
	\end{align} leading to (i). Similarly $E_{-}>E_{+}$ leads to (ii). Thus, assuming \eqref{eq:edge1} we have
	\beq
	M_{\rho}(E_{0})=1-\left(\int_{\R_{+}}\frac{x}{x-\Omega_{2}(E_{0})}\dd\nu(x)\right)^{-1}<1,
	\eeq
	yet at the same time
	\beq
	M_{\rho}(E_{0})=1-\left(\int_{\R_{+}}\frac{x}{x-\Omega_{1}(E_{0})}\dd\mu(x)\right)^{-1}>1,
	\eeq
	which is a contradiction.
	Also for the other case in (iii) the same argument gives a contradiction, concluding the proof.
\end{proof}

\begin{remark}\label{rem:single0}
	Even if we allow the measures to touch $0$, we have the same characterization \eqref{eq:edgechar} for $E\in\caV\setminus\{0\}$. Also in this case, we still have exactly one (no, respectively) solution in case (ii) ((iii), respectively) in the proof of Proposition \ref{prop:edgechar2}. Using a similar argument as in Remark \ref{rem:sqrt0} we get $E_{-}=0$, so that $\{x\in\R:f(x)>0\}=(0,E_{+})$ is still a single interval.
\end{remark}

As a final ingredient, we prove that the subordination functions $\Omega_{2}$ and $\Omega_{1}$ have square root behavior at the edges $E_{-}$ and $E_{+}$, so that $M_{\rho}=M_{\nu}\circ\Omega_{2}=M_{\mu}\circ\Omega_{1}$ also does.
\begin{proposition}\label{prop:sqrt}
	Let $\mu$ and $\nu$ satisfy the assumptions of Theorem \ref{thm:main} and let $\supp\rho=[E_{-},E_{+}]$. Then there exist strictly positive constants $\gamma_{1,-},\gamma_{1,+},\gamma_{2,-}$, and $\gamma_{2,+}$ such that for $j=1,2$, 
	\beq\label{eq:sqrtsuborl}
	\Omega_{j}(z)=\Omega_{j}(E_{-})+\gamma_{j,-}\sqrt{E_{-}-z}+O(\absv{z-E_{-}}^{3/2}),
	\eeq
	for $z$ in a neighborhood of $E_{-}$ with the principal branch of square root(with $\sqrt{-1}=\ii$). Similarly for $E_{+}$, we have
	\beq\label{eq:sqrtsuboru}
	\Omega_{j}(z)=\Omega_{j}(E_{+})+\gamma_{j,+}\sqrt{z-E_{+}}+O(\absv{z-E_{-}}^{3/2}),
	\eeq
	for $z$ in a neighborhood of $E_{+}$ with the same branch of square root.
\end{proposition}
\begin{proof}
	For simplicity, we focus on the behavior of $\Omega_{2}$. The corresponding results for $\Omega_{1}$ can be proved analogously. We first note that
	\beq
	M_{\mu}'(\Omega)=\frac{1}{(\Omega m_{\mu}(\Omega)+1)^{2}}\int_{\R_{+}}\frac{x}{(x-\Omega)^{2}}\dd\mu\neq 0,
	\eeq
	whenever $\Omega\notin\supp\mu$. Since  $\Omega_{1}(E_{-})\in (0,E_{-}^{\mu})$, the inverse $M_{\mu}^{-1}$ is well-defined and analytic in a neighborhood of $M_{\mu}\circ\Omega_{1}(E_{-})=M_{\nu}\circ\Omega_{2}(E_{-})$, hence the function
	\beq
	\wt{z}_{-}(\Omega)= \Omega\frac{M_{\mu}^{-1}\circ M_{\nu}(\Omega)}{M_{\nu}(\Omega)}.
	\eeq
	is analytic in a neighborhood of $\Omega_{2}(E_{-})$. Furthermore we find that if $z$ is sufficiently close to $E_{-}$,
	\beq
	M_{\mu}^{-1}\circ M_{\nu}\circ \Omega_{2}(z)=M_{\mu}^{-1}\circ M_{\mu}\circ\Omega_{1}(z)=\Omega_{1}(z),
	\eeq so that
	\beq
	\wt{z}_{-}(\Omega_{2}(z))
	=\Omega_{2}(z)\frac{M_{\mu}^{-1}\circ M_{\nu}\circ\Omega_{2}(z)}{M_{\rho}(z)} 
	=\frac{\Omega_{2}(z)\Omega_{1}(z)}{M_{\rho}(z)}=z.
	\eeq 
	
	Now given the fact that $\wt{z}_{-}$ is analytic in a neighborhood of $\Omega_{2}(E_{-})$, we consider its Taylor expansion 
	\beq
	\wt{z}_{-}(\Omega)=E_{-}+\wt{z}_{-}'(\Omega_{2}(E_{-}))(\Omega-\Omega_{2}(E_{-}))+\frac{1}{2}\wt{z}_{-}''(\Omega_{2}(E_{-}))(\Omega-\Omega_{2}(E_{-}))^{2}+R_{-}(\Omega)
	\eeq
	around $\Omega_{2}(E_{-})$, where $R_{-}(\Omega)=O(\absv{\Omega-\Omega_{2}(E_{-})}^{3})$. We shall prove that 
	\begin{align}
	&\wt{z}_{-}'(\Omega_{2}(E_{-}))=0, &&\wt{z}_{-}''(\Omega_{2}(E_{-}))\neq 0.
	\end{align}
	Along the remaining proof, values at $\Omega_{2}(E_{\pm})$ of the function $M_{\nu}$ and its derivatives are used repeatedly. They can be derived from Proposition \ref{prop:stab} combined with either the definition of $M_{\mu}$ or Lemmas~\ref{lem:NPrepM} and \ref{lem:NPrepMsupp}. For readers' convenience, they are listed below:
	\begin{align}\label{eq:dervalues}
	&0<\Omega_{1}(E_{-})<E_{-}^{\mu}, &&E_{+}^{\mu}<\Omega_{1}(E_{+}), &&0<\Omega_{2}(E_{-})<E_{-}^{\nu}, &&E_{+}^{\nu}<\Omega_{2}(E_{+});\nonumber\\
	&M_{\nu}(\Omega_{2}(E_{-}))\in (0,1),	&&M_{\nu}(\Omega_{2}(E_{+}))\in (1,\infty),	&&M_{\nu}'(\Omega_{2}(E_{-}))>1, &&M_{\nu}'(\Omega_{2}(E_{+}))>1; \nonumber\\
	&H_{\nu}'(\Omega_{2}(E_{-}))>0,	&&H_{\nu}'(\Omega_{2}(E_{+}))>0,	&&M_{\nu}''(\Omega_{2}(E_{-}))>0,	&&M_{\nu}''(\Omega_{2}(E_{+}))<0.
	\end{align}
	For the first derivative of $\wt{z}_{-}$, by the definition of $\wt{z}_{-}$ and the fact that $E_{-}\in\caV$ satisfies \eqref{eq:edgechar}, we see that
	\begin{multline}\label{eq:1stder=0}
	\wt{z}_{-}'(\Omega_{2}(E_{-}))
	=\left[\frac{M_{\mu}^{-1}\circ M_{\nu}(\Omega)}{M_{\nu}(\Omega)}
	+\Omega\frac{M_{\nu}'(\Omega)}{M_{\nu}(\Omega)M'_{\mu}(M_{\mu}^{-1}\circ M_{\nu})(\Omega)}
	-\Omega\frac{(M_{\mu}^{-1}\circ M_{\nu})(\Omega)}{M_{\nu}(\Omega)^{2}}M_{\nu}'(\Omega)\right]_{\Omega=\Omega_{2}(E_{-})}	\\
	=\frac{1}{M_{\mu}'(\Omega_{1}(E_{-}))}
	\Big(\frac{\Omega_{1}(E_{-})}{M_{\mu}(\Omega_{1}(E_{-}))}M_{\mu}'(\Omega_{1}(E_{-}))
	+\frac{\Omega_{2}(E_{-})}{M_{\nu}(\Omega_{2}(E_{-}))}M_{\nu}'(\Omega_{2}(E_{-}))\\
	-\frac{\Omega_{2}(E_{-})\Omega_{1}(E_{-})}{M_{\rho}(E_{-})^{2}}M_{\mu}'(\Omega_{1}(E_{-}))M_{\nu}'(\Omega_{2}(E_{-}))\Big)=0,
	\end{multline}
	since $M_{\mu}'(\Omega_{1}(E_{-}))\neq 0$. 
	
	Differentiating once again we obtain
	\begin{multline}\label{eq:2ndder0}
	\wt{z}_{-}''(\Omega)=2\left[\frac{M_{\nu}'(\Omega)}{M_{\nu}(\Omega)M'_{\mu}(M_{\mu}^{-1}\circ M_{\nu})(\Omega)}
	-\frac{(M_{\mu}^{-1}\circ M_{\nu})(\Omega)}{M_{\nu}(\Omega)^{2}}M_{\nu}'(\Omega)
	-\Omega\frac{M_{\nu}'(\Omega)^{2}}{M_{\nu}(\Omega)^{2}(M_{\mu}'\circ M_{\mu}\circ M_{\nu})(\Omega)}\right] \\
	+\frac{\Omega}{M_{\nu}(\Omega)}
	\left[\frac{M_{\nu}''(\Omega)}{(M_{\mu}'\circ M_{\mu}^{-1}\circ M_{\nu})(\Omega)}
	-\frac{M_{\nu}'(\Omega)(M_{\mu}''\circ M_{\mu}^{-1}\circ M_{\nu})(\Omega)}{(M_{\mu}'\circ M_{\mu}^{-1}\circ M_{\nu})(\Omega)^{3}}\right] \\
	-\Omega (M_{\mu}^{-1}\circ M_{\nu})(\Omega)\left[\frac{M_{\nu}''(\Omega)}{M_{\nu}(\Omega)^{2}}-2\frac{M_{\nu}'(\Omega)^{2}}{M_{\nu}(\Omega)^{3}}\right].
	\end{multline}
	Plugging in $\Omega=\Omega_{2}(z)$, the second derivative in \eqref{eq:2ndder0} simplifies to
	\begin{multline}\label{eq:2ndder1}
	\wt{z}_{-}''(\Omega_{2}(z))=2\left[\frac{M_{\nu}'(\Omega_{2}(z))}{M_{\rho}(z)M_{\mu}'(\Omega_{1}(z))}
	-\frac{\Omega_{1}(z)M_{\nu}'(\Omega_{2}(z))}{M_{\rho}(z)^{2}}
	-\frac{\Omega_{2}(z)M_{\nu}'(\Omega_{2}(z))^{2}}{M_{\rho}(z)^{2}M_{\mu}'(\Omega_{1}(z))}\right] \\
	+\frac{\Omega_{2}(z)}{M_{\rho}(z)}
	\left[\frac{M_{\nu}''(\Omega_{2}(z))}{M_{\mu}'(\Omega_{1}(z))}
	-\frac{M_{\nu}'(\Omega_{2}(z))^{2}M_{\mu}''(\Omega_{1}(z))}{M_{\mu}'(\Omega_{1}(z))^{3}}\right]
	-\Omega_{2}(z)\Omega_{1}(z)\left[\frac{M_{\nu}''(\Omega_{2}(z))}{M_{\rho}(z)^{2}}-2\frac{M_{\nu}'(\Omega_{2}(z))^{2}}{M_{\rho}(z)^{3}}\right].
	\end{multline}
	Combining the first three and the last terms on the right-hand side of \eqref{eq:2ndder1} we have
	\begin{multline}
	\frac{2M_{\nu}'(\Omega_{2}(z))}{M_{\rho}(z)M_{\mu}'(\Omega_{1}(z))}\left(1-\frac{\Omega_{1}(z)}{M_{\rho}(z)}M_{\mu}'(\Omega_{1}(z))-\frac{\Omega_{2}(z)}{M_{\rho}(z)}M_{\nu}'(\Omega_{2}(E_{-}))+\frac{\Omega_{1}(z)\Omega_{2}(z)}{M_{\rho}(z)^{2}}M_{\mu}'(\Omega_{1}(z))M_{\nu}'(\Omega_{2}(z))\right)\\
	=\frac{2M_{\nu}'(\Omega_{2}(z))}{M_{\rho}(z)^{3}M_{\mu}'(\Omega_{1}(z))}\times\left(\text{left-hand side of \eqref{eq:edgechar}}\right),
	\end{multline}
	which is zero at $z=E_{-}$ by $\eqref{eq:edgechar}$. On the other hand, collecting the fourth and sixth terms on the right-hand side of \eqref{eq:2ndder1} we have
	\begin{multline}
	\frac{\Omega_{2}(z)M_{\nu}''(\Omega_{2}(z))}{M_{\rho}(z)M_{\mu}'(\Omega_{1}(z))}-\frac{\Omega_{2}(z)\Omega_{1}(z)M_{\nu}''(\Omega_{2}(z))}{M_{\rho}(z)^{2}}\\	=-\frac{\Omega_{1}(z)^{2}\Omega_{2}(z)M_{\nu}''(\Omega_{2}(z))}{M_{\rho}(z)^{2}M_{\mu}'(\Omega_{2}(z))}\left(\frac{M_{\mu}'(\Omega_{1}(z))}{\Omega_{1}(z)}-\frac{M_{\rho}(z)}{\Omega_{1}(z)^{2}}\right)
	=-\frac{\Omega_{1}(z)^{2}\Omega_{2}(z)M_{\nu}''(\Omega_{2}(z))}{M_{\rho}(z)^{2}M_{\mu}'(\Omega_{2}(z))}H_{\mu}'(\Omega_{1}(z)),
	\end{multline}
	which is negative at $z=E_{-}$ by \eqref{eq:dervalues}.	Therefore, we have
	\beq\label{eq:2ndder}
	\wt{z}_{-}''(\Omega_{2}(E_{-}))
	=-\frac{\Omega_{1}(E_{-})^{2}\Omega_{2}(E_{-})M_{\nu}''(\Omega_{2}(E_{-}))}{M_{\rho}(E_{-})^{2}M_{\mu}'(\Omega_{2}(E_{-}))}H_{\mu}'(\Omega_{1}(E_{-}))
	-\frac{\Omega_{2}(E_{-})M_{\nu}'(\Omega_{2}(E_{-}))^{2}M_{\mu}''(\Omega_{1}(E_{-}))}{M_{\rho}(E_{-})M_{\mu}'(\Omega_{1}(E_{-}))^{3}}<0.
	\eeq
	
	Defining $\gamma_{2,-}=\sqrt{-2/(\wt{z}_{-}''(\Omega_{2}(E_{-})))}>0$, from the Taylor series expansion of $\wt{z}_{-}(\Omega)$ we see that
	\beq
	z-E_{-}=-\frac{1}{\gamma_{2,-}^{2}}(\Omega_{2}(z)-\Omega_{2}(E_{-}))^{2}+R_{-}(\Omega_{2}(z))
	\eeq
	for $z$ in a neighborhood of $E_{-}$,  by continuity of $\Omega_{2}$. Inverting the expansion concludes the proof for lower edge. 
	
	For the upper edge, we similarly define $\wt{z}_{+}$ to be the inverse function of $\Omega_{2}(z)$ in a neighborhood of $E_{+}$, so that its derivatives have exactly the same forms as those of $\wt{z}_{-}$. Thus $\wt{z}_{+}'(\Omega_{2}(E_{+}))=0$ can be proved in a completely analogous manner. On the other hand, observing that $M_{\nu}''(\Omega_{2}(E_{+}))$ and $M_{\mu}''(\Omega_{1}(E_{+}))$ are strictly negative, we can immediately see that $\wt{z}_{+}''(\Omega_{2}(E_{+}))>0$ from \eqref{eq:2ndder}. Now the result follows again from Taylor expansion of $\wt{z}_{+}$ around $\Omega_{2}(E_{+})$. Note that the difference of signs of $\wt{z}''_{\pm}$ at $\Omega_{2}(E_{\pm})$ induces that of branches of square roots in \eqref{eq:sqrtsuborl} and \eqref{eq:sqrtsuboru}.
\end{proof}

\begin{proof}[Proof of Theorem~\ref{thm:main}]
	We first recall that,
	\beq
	\im (zm_{\rho}(z)+1)=\im (\Omega_{2}(z)m_{\nu}(z)+1)= I_{\nu}(\Omega_{2}(z))\im \Omega_{2}(z), \quad \forall z\in\C_{+}.
	\eeq
	From Proposition~\ref{prop:stab}, we find that $I_{\nu}(\Omega_{2}(z))$ is bounded below and above uniformly in any compact set $\caD\subset\C_{+}\cup\R\setminus\{0\}$. Thus we see that
	\beq
	\frac{xf(x)}{\sqrt{x-E_{-}}}=\frac{\im (xm_{\rho}(x)+1)}{\sqrt{x-E_{-}}}
	\eeq
	is bounded above and below by
	\beq
	\frac{\im \Omega_{2}(x)}{\sqrt{x-E_{-}}}	
	=
	\begin{cases}
		\gamma_{2,-}\frac{\im \sqrt{E_{-}-x}}{\sqrt{x-E_{-}}}+O(\absv{z-E_{-}}) =\gamma_{2,-}+O(\absv{z-E_{-}}) & \text{if }x>E_{-},\\ 
		0 & \text{otherwise},
	\end{cases}
	\eeq 
	for $x$ around $E_{-}$. Similar reasoning for $E_{+}$ proves the assertion.
\end{proof}
\begin{remark}\label{rem:sqrt0pf}
	As easily seen, the proof Theorem \ref{thm:main} for the upper and lower edge are totally separated. Thus even if $E_{-}^{\mu}$ or $E_{-}^{\nu}$ is $0$, since all the other ingredient for the proof is intact around $E_{+}$, we still get the same result at the upper edge.
\end{remark}

\subsubsection*{Acknowledgments}
The author would like to thank Ji Oon Lee for helpful suggestions and discussion. The author is grateful to the anonymous referees for carefully reading our manuscript and providing valuable comments. This research has been supported by TJ Park Science Fellowship of POSCO TJ Park Foundation.


\begin{thebibliography}{10}
	
	\bibitem{Akhiezer1965}
	N.~I. Akhiezer.
	\newblock {\em The classical moment problem and some related questions in
		analysis}.
	\newblock Translated by N. Kemmer. Hafner Publishing Co., New York, 1965.
	
	\bibitem{Bao-Erdos-Schnelli2016}
	Z.~Bao, L.~Erd\H{o}s, and K.~Schnelli.
	\newblock Local stability of the free additive convolution.
	\newblock {\em J. Funct. Anal.}, 271(3):672--719, 2016.
	
	\bibitem{Bao-Erdos-Schnelli2017arXiv}
	Z.~{Bao}, L.~{Erd\H{o}s}, and K.~{Schnelli}.
	\newblock {Spectral rigidity for addition of random matrices at the regular
		edge}.
	\newblock {\em arXiv e-prints}, page arXiv:1708.01597, Aug 2017.
	
	\bibitem{Bao-Erdos-Schnelli2018}
	Z.~{Bao}, L.~{Erd\H{o}s}, and K.~{Schnelli}.
	\newblock {On the support of the free additive convolution}.
	\newblock {\em arXiv e-prints}, page arXiv:1804.11199, Apr 2018.
	
	\bibitem{Belinschi2003}
	S.~T. Belinschi.
	\newblock The atoms of the free multiplicative convolution of two probability
	distributions.
	\newblock {\em Integral Equations Operator Theory}, 46(4):377--386, 2003.
	
	\bibitem{Belinschi2005}
	S.~T. Belinschi.
	\newblock {\em Complex analysis methods in noncommutative probability}.
	\newblock ProQuest LLC, Ann Arbor, MI, 2005.
	\newblock Thesis (Ph.D.)--Indiana University.
	
	\bibitem{Belinschi2006}
	S.~T. Belinschi.
	\newblock A note on regularity for free convolutions.
	\newblock {\em Ann. Inst. Henri Poincar\'{e} Probab. Stat.}, 42(5):635--648,
	2006.
	
	\bibitem{Belinschi2008}
	S.~T. Belinschi.
	\newblock The {L}ebesgue decomposition of the free additive convolution of two
	probability distributions.
	\newblock {\em Probab. Theory Related Fields}, 142(1-2):125--150, 2008.
	
	\bibitem{Belinschi2014}
	S.~T. Belinschi.
	\newblock {$L^\infty$}-boundedness of density for free additive convolutions.
	\newblock {\em Rev. Roumaine Math. Pures Appl.}, 59(2):173--184, 2014.
	
	\bibitem{Belinschi-Bercovici2007}
	S.~T. Belinschi and H.~Bercovici.
	\newblock A new approach to subordination results in free probability.
	\newblock {\em J. Anal. Math.}, 101:357--365, 2007.
	
	\bibitem{Biane1998}
	P.~Biane.
	\newblock Processes with free increments.
	\newblock {\em Math. Z.}, 227(1):143--174, 1998.
	
	\bibitem{Bingham-Goldie-Teugels1989}
	N.~H. Bingham, C.~M. Goldie, and J.~L. Teugels.
	\newblock {\em Regular variation}, volume~27 of {\em Encyclopedia of
		Mathematics and its Applications}.
	\newblock Cambridge University Press, Cambridge, 1989.
	
	\bibitem{Lee-Schnelli2013}
	J.~O. Lee and K.~Schnelli.
	\newblock Local deformed semicircle law and complete delocalization for
	{W}igner matrices with random potential.
	\newblock {\em J. Math. Phys.}, 54(10):103504, 62, 2013.
	
	\bibitem{HBMF2010}
	F.~W.~J. Olver, D.~W. Lozier, R.~F. Boisvert, and C.~W. Clark, editors.
	\newblock {\em N{IST} handbook of mathematical functions}.
	\newblock U.S. Department of Commerce, National Institute of Standards and
	Technology, Washington, DC; Cambridge University Press, Cambridge, 2010.
	\newblock With 1 CD-ROM (Windows, Macintosh and UNIX).
	
	\bibitem{Voiculescu1986}
	D.~Voiculescu.
	\newblock Addition of certain noncommuting random variables.
	\newblock {\em J. Funct. Anal.}, 66(3):323--346, 1986.
	
	\bibitem{Voiculescu1987}
	D.~Voiculescu.
	\newblock Multiplication of certain noncommuting random variables.
	\newblock {\em J. Operator Theory}, 18(2):223--235, 1987.
	
	\bibitem{Voiculescu1991}
	D.~Voiculescu.
	\newblock Limit laws for random matrices and free products.
	\newblock {\em Invent. Math.}, 104(1):201--220, 1991.
	
	\bibitem{Voiculescu1993}
	D.~Voiculescu.
	\newblock The analogues of entropy and of {F}isher's information measure in
	free probability theory. {I}.
	\newblock {\em Comm. Math. Phys.}, 155(1):71--92, 1993.
	
\end{thebibliography}
\end{document}